\documentclass[12pt,reqno]{amsart}

\usepackage{amsmath,amsfonts,amssymb,amsthm,mathtools}
\usepackage{hyperref}
\usepackage{gensymb}
\usepackage{graphicx}	
\usepackage{float}
\usepackage{color}
\usepackage{longtable}

\usepackage{pdfpages}
\usepackage{footnote}
\usepackage{appendix}
\usepackage{sidecap}
\usepackage{tabularx}
\usepackage{framed}
\makesavenoteenv{tabular}
\makesavenoteenv{table}
\usepackage{array}
\usepackage{booktabs}

\usepackage{tikz}
\usepackage{tikz-cd}
\usetikzlibrary{shapes.geometric,calc,positioning,3d,intersections,arrows}
\usetikzlibrary{decorations.markings}
\usepackage{pgfplots}

\calclayout

\renewenvironment{framed}[1][\hsize]
   {\MakeFramed{\hsize#1\advance\hsize-\width \FrameRestore}}%
   {\endMakeFramed}

\newcommand{\la}{\langle}
\newcommand{\ra}{\rangle}

\renewcommand\a{\alpha}

\newcommand\bx{{\boldsymbol x}}
\newcommand\by{{\boldsymbol y}}
\newcommand\bz{{\boldsymbol z}}
\newcommand\bu{{\boldsymbol u}}

\newcommand\g{\mathfrak{g}}
\newcommand\gl{\mathfrak{gl}}

\newcommand\oo{\omega}
\newcommand\op[1]{\mathop{\rm #1}\nolimits}
\newcommand\p{\partial}

\newcommand\R{{\mathbb R}}
\renewcommand\sl{\mathfrak{sl}}
\renewcommand\sp{\mathfrak{sp}}
\newcommand\we{\wedge}

\newtheorem{thm}{Theorem}[section] 
\theoremstyle{definition}

\newtheorem{theorem}[thm]{Theorem}

\pgfplotsset{compat=1.15}

\begin{document}
\title{Differential Invariants of Linear Symplectic Actions}
\author{J\o rn Olav Jensen}
\author{Boris Kruglikov}
 \address{Institute of Mathematics and Statistics, UiT the Arctic University of Norway, Troms\o\ 90-37, Norway.
E-mails: \textsc{jje041@post.uit.no,\ boris.kruglikov@uit.no}.}
\keywords{Differential invariants, invariant derivations, symplectic and contact spaces}
\subjclass[2010]{53A55, 14L24; 37J15, 15A72}
 \maketitle

\begin{abstract}
We consider the equivalence problem for symplectic and conformal symplectic group actions on submanifolds and functions of
symplectic and contact linear spaces. This is solved by computing differential invariants via the Lie-Tresse theorem.
\end{abstract}

 \section{Introduction}\label{S1}

Differential invariants of various groups play an important role in applications \cite{Ov,Ol,T}.
Classical curvatures of submanifolds in Euclidean space arise as
differential invariants of the orthogonal group.
The corresponding problem for symplectic spaces was initiated in \cite{CW}.
Further works in this direction include \cite{D,KOT,MK,MH,MN,V}.
In this paper we consider the linear symplectic group action and compute the corresponding algebra of differential invariants.
We will use either the standard representation or its trivial extension; other actions were also considered in the
literature \cite{X} and we comment on the relations of the above cited works to ours at the conclusion of the paper.

Let $V=\R^{2n}(\bx,\by)$ be equipped with the standard symplectic form $\omega=\sum_1^n dx_i\wedge dy_i$.
Every infinitesimal symplectic transformation of $V$ is given by the Hamiltonian function $H\in C^\infty(V)$ and has the form
$X_H=\oo^{-1}dH$, and the Lie bracket of vector fields corresponds to the Poisson bracket of functions.
By the Darboux-Givental theorem the action of $\op{Symp}(V,\oo)$ has no local invariants. However
these arise when we restrict to finite-dimensional subalgebras/subgroups. Namely, functions $H$ quadratic in $x,y$ form a subalgebra
isomorphic to $\sp(2n,\R)$. For functions of degree $\leq2$ we get the affine symplectic algebra $\sp(2n,\R)\ltimes\R^{2n}$.
We will concentrate on the linear case and compute the algebra of differential invariants for submanifolds and functions on $V$.

It turns out that for curves and hypersurfaces one can describe the generators for all $n$ that we provide,
while for the case of dimension and codimension greater than one, this becomes more complicated. Of those we consider in details
only the case of surfaces in $\R^4$.
Generators of the  algebra of differential invariants will be presented in the Lie-Tresse form as functions and derivations, and for lower
dimensions we also compute the differential syzygies. We will mainly discuss the geometric coordinate-free approach. The explicit formulae are rather
large and will be shown in the appendix only for $n=2$.

We also consider the space $W=\R^{2n+1}(\bx,\by,z)$ equipped with the standard contact form $\a=dz-\sum_1^n y_idx_i$.
Every infinitesimal contact transformation of $W$ is given by the contact Hamiltonian $H\in C^\infty(W)$ via
$\a(X_H)=H$, $X_H(\a)=\p_u(H)$, and the Lie bracket of vector fields corresponds to the Lagrange bracket of functions.
Again, the action of $\op{Cont}(W,[\a])$ has no local invariants, however
these arise when we restrict to finite-dimensional subalgebras/subgroups. Namely, functions $H$ quadratic in $\bx,\by,z$
with weights $w(x_i)=1=w(y_i)$, $w(z)=2$ form a subalgebra isomorphic to $\mathfrak{csp}(2n,\R)$. For functions of degree $\leq2$
we get the affine extension $(\R\oplus\sp(2n,\R))\ltimes\R^{2n+1}$.
For simplicity, we will concentrate on the action of $\sp(2n,\R)$, and then comment how to extend to the conformally symplectic algebra
and include the translations.

It is interesting to remark on the computational aspect of the results. There are two approaches to compute differential invariants.
The infinitesimal method is based on the defining Lie equations and
works universally for arbitrary Lie algebras of vector fields.
The moving frame method is based on elimination of group parameters and
is dependent on explicit parametrization
of the Lie group (or pseudogroup in infinite-dimensional situation) and its action. In \textsc{Maple} these in turn rely on \verb"pdsolve"
and \verb"eliminate" commands or some algorithmically optimized versions of those (via Gr\"obner basis or similar).
For the problem at hand we can use both since one can locally
parametrize the group $\op{Sp}(2n,\R)$ and its linear action.
The Lie algebra method works well in dimension 2 (symplectic case $n=1$) and fails further.
The Lie group method works well in dimension 3 (contact case $n=1$) and fails further.
Computational difficulties obstruct finishing calculations already in dimension 4 with these straightforward approaches.
We show however how other geometric methods allow to proceed further.

This paper is partially based on the results of \cite{JJ},
extending and elaborating it in several respects.
Some applications will be briefly discussed at the end of the paper. The paper is organized as follows.
In the next section we recall the basics. Then we describe in turn differential invariants of functions, curves and hypersurfaces
in symplectic vector spaces, and also discuss the particular case of surfaces in $\R^4$. Then we briefly discuss the invariants in contact vector spaces
and demonstrate how to compute differential invariants for conformal and affine extensions from our preceding computations.

We present most computations explicitly.
Some large formulae are delegated to the appendix, the other can be found (as maple and pdf) in ancillary files of this arXiv.

\section{Recollections and setup}

We refer to \cite{KL0} for details of the jet-formalism,
summarizing the essentials here.

\subsection{Jets}

Let $M$ be a smooth manifold. Two germs at $a\in M$ of
submanifolds $N_1,N_2\subseteq M$ of dimension $n$ and codimension $m$
are equivalent if they are tangent up to order $k$ at $a$.
The equivalence class $[N]_a^k$ is called the $k$-jet of $N$ at $a$.
Denote $J_a^k(M,n)$ the set of all $k$-jets at $a$ and
$J^k(M,n)=\cup_{a\in M}J_a^k(M,n)$ the space of $k$-jets of $n$-submanifolds.
This is a smooth manifold of dimension $n+m\binom{n+k}{k}$
and there are natural bundle projections $\pi_{k,l}:J^k(M,n)\to J^l(M,n)$
for $k>l\ge 0$. Note that $J^0(M,n)=M$ and $J^1(M,n)=\op{Gr}_n(TM)$,
while $\pi_{k,k-1}:J^k(M,n)\to J^{k-1}(M,n)$ are affine bundles for $k>1$.

Since functions $f\in C^\infty(M)$ can be identified with their graphs
$\Sigma_f\subset M\times\R$, the space of $k$-jets of functions $J^kM$
is defined as the space of $k$-jets of hypersurfaces $\Sigma\subset M\times\R$
transversal to the fibers of the projection to $M$.
This jet space embeds as an open subset into $J^k(M\times\R,n)$, where $n=\dim M$
(and $m=1$) and so its dimension is $n+\binom{n+k}{k}$.

Sometimes we denote spaces $J^kM$ and $J^k(M,n)$ simply by $J^k$.
The inverse limit along projections $\pi_{k,k-1}$ yields the space
$J^\infty=\varprojlim J^k$.

In local coordinates $(\bx,\by)$ on $M$ a submanifold $N$ can be written as
$y^j=y^j(x^i)$, $i=1,\dots,n$, $j=1,\dots,m$. Then the jet-coordinates
are given by $x^i([N]_a^k)=a^i$, $y^j_\sigma([N]_a^k)=
\frac{\partial^{|\sigma|}y^j}{\partial x^\sigma}(a)$
for a multi-index $\sigma=(i_1,\hdots,i_n)$ of length
$|\sigma|=\sum_1^n i_s\le k$.

For the jets of functions $u=u(\bx)$ we use the jet-coordinates
$x^i([u]_a^k)=a^i$,
$u_{\sigma}([u]_a^k)=\frac{\p^{|\sigma|}u}{\p x^\sigma}(a)$.
We sometimes also write $u$ instead of $u_0$, and we often lower indices
for the base coordinates, like $x_i$ instead of $x^i$ etc,
if no summation suffers.

\subsection{Prolongations}

A Lie group action on a manifold $M$ is a homomorphism
$\Phi:G\to\op{Diff}(M)$. Any $g\in G$ determines a point transformation
$\Phi_g(a)=g\cdot a$, $a\in M$. This induces an action on germs of submanifolds,
hence on jets of submanifolds, namely
 $$
\Phi_g^{(k)}([N]_a^k)=[\Phi_g(N)]_{\Phi_g(a)}^k.
 $$

Similarly, if $X$ is a vector field on $M$, corresponding to the
Lie algebra $\g=\op{Lie}(G)$, the prolongation
or lift gives a vector field $X^{(k)}$ on $J^k$.
If $(\bx,\bu)$ are local coordinates on $M$
(with $x^i$ interpreted as independent and $u^j$ as dependent variables)
and the vector field is given as $X=a^i\partial_{x^i}+b^j\partial_{u^j}$,
then its prolongation has the form
 $$
X^{(k)} = a^i\mathcal{D}_{x^i}^{(k+1)}+\sum_{|\sigma|\le k} \mathcal{D}_{\sigma}(\varphi^j)\partial_{u_{\sigma}^j},
 $$
where $\varphi=(\varphi^1,\hdots,\varphi^m)$ and $\varphi^j=b^j-a^iu_i^j$
is the generating vector-function, $\mathcal{D}_{x^i}=\p_{x^i}+\sum_{j,\tau}u_{\tau+1_i}^j\p_{u_\tau^j}$
is the total derivative, $\mathcal{D}_{x^i}^{(k+1)}$ its truncation
(restriction to $(k+1)$-jets: $|\tau|\leq k$) and $\mathcal{D}_{\sigma}=\mathcal{D}_{x^1}^{i_1}\cdots\mathcal{D}_{x^n}^{i_n}$
for $\sigma=(i_1,\dots,i_n)$ is the iterated total derivative.

\subsection{Differential Invariants}

A differential invariant of order $k$ is a function $I$ on $J^k$, which
is constant on the orbits of $\Phi^{(k)}$ action of $G$.
If the Lie group $G$ is connected this is equivalent to
$L_{X^{(k)}}I=0$ for all $X\in\g$ (some care should be taken with this
statement, mostly related to usage of local coordinate charts in jets,
see the first example in \cite{KL}).

The space of $k$-th order differential invariants forms a commutative
algebra over $\mathbb{R}$, denoted by $\mathcal{A}_k$.
The injection $\pi_{k+1,k}^*$ induces the embedding $\mathcal{A}_k\subset\mathcal{A}_{k+1}$, and in the inductive limit
we get the algebra of differential invariants
$\mathcal{A}\subseteq C^{\infty}(J^{\infty})$, namely
 $$
\mathcal{A}=\varinjlim\mathcal{A}_{k}.
 $$

Denote by $G_a=\{g\in G:g\cdot a=a\}$ the stabilizer of $a\in M$.
This subgroup of $G$ acts on $J^k_a$. The prolonged action of $G$
is called algebraic if the prolongation $G_a^{(k)}$ is an algebraic group
acting algebraically on $J^k_a$ $\forall$ $a\in M$.
For our problem the action of $G$ on $M$ is almost transitive and
algebraic, so by \cite{KL} the invariants $I\in\mathcal{A}$
can be taken as {\em rational functions} in jet-variables $u^j_\sigma$;
moreover they may be chosen {\bf polynomial} starting from some jet-order.
This will be assumed in what follows.

In our situation $\mathcal{A}$ is not finitely generated in the usual sense since the number of independent invariants is infinite.
We will use the Lie-Tresse theorem \cite{KL} that guarantees that
$\mathcal{A}$ is generated by a finite set of differential invariants
and invariant derivations.

Recall that an invariant derivation is such a horizontal (or Cartan)
derivation $\nabla:\mathcal{A}\to\mathcal{A}$ (obtained by a combination of
total derivatives) that it commutes with the action of the group:
$\forall g\in G$ we have $g_{*}^{(k+1)}\nabla=\nabla g_{*}^{(k)}$
for $k\ge k_0$, where $k_0$ is the order of $\nabla$, which can be identified
with the highest order of coefficients in the decomposition
$\nabla=\sum_ia^i(x,u^j_\sigma)\mathcal{D}_{x^i}$. Equivalently we can write
$\forall X\in\g$: $\mathcal{L}_{X^{(k+1)}}\nabla=\nabla\mathcal{L}_{X^{(k)}}$ for $k\ge k_0$. This implies $\nabla:\mathcal{A}_k\to\mathcal{A}_{k+1}$
in the same range.

Invariant derivations form a submodule
$\mathcal{CD}^G\subseteq\mathcal{CD}(J^\infty)$ in the space of
all horizontal derivations. It is a finitely
generated $\mathcal{A}$ module: any $\nabla\in\mathcal{CD}^G$ has the form
$\nabla=I^i\nabla_i$ for a fixed set $\nabla_i$ and $I^i\in\mathcal{A}$.
By \cite[Theorem 21]{KL} the number of derivations $\nabla_i$ is $n$.

We compose iterated operators $\nabla_{J}:\mathcal{A}_k\to\mathcal{A}_{k+|J|}$
for multi-indices $J$, and then $\mathcal{A}$ is generated by
$\nabla_{J}I_i$ for a finite set of $I_i$.

\subsection{Counting the invariants}

An important part of our computations is a count of independent
differential invariants. Denote the number of those on the level
of $k$-jets by $s_k$. This number is equal to the transcendence degree
of the field of differential invariants on $J^k$
(when the elements of $\mathcal{A}_k$ are rational functions) and it coincides with the codimension of $G^{(k)}$ orbit in $J^k$.

Since in our case $G$ is a (finite-dimensional) Lie group, the action becomes eventually free, i.e.\ $G^{(k)}_{a_k}=\op{id}$ for
sufficiently large $k\ge\ell$ and generic $a_k\in J^k$ cf.\ \cite{Ol}. In this case the orbit is diffeomorphic to $G$, in particular
$s_k=\dim J^k-\dim G$ for $k\ge\ell$.

The number of "pure order" $k$ differential invariants is $h_k=s_k-s_{k-1}$, so it becomes 
 $$
h_k=\dim J^k-\dim J^{k-1}=m\tbinom{n+k-1}k\ \text{ for }\ k>\ell.
 $$
The Poincar\'e function $P(z)=\sum_{k=0}^\infty h_kz^k$ is rational in all local problems of analysis according to Arnold's conjecture \cite{Kr}.
In our case this $P(z)$ differs from $m(1-z)^{-n}$ by a polynomial reflecting the action of $G$.

Note that by the eventual freeness of the action, the algebra $\mathcal{A}$ is generated by invariants and derivations at most from the jet-level
$\ell$.

\subsection{The equivalence problem}

The generators $I_i$ ($1\le i\le s$), $\nabla_j$ ($1\le j\le n$) are not independent, i.e.\ the algebra $\mathcal{A}$ is not freely generated by them, in general.
A differential syzygy is a relation among these generators. Such an expression has the form $F(\nabla_{J_1}(I_{i_1}),\hdots,\nabla_{J_n}(I_{i_t}))=0$,
where $F$ is a function of $t$ arguments and $J_1,\hdots,J_t$ are multi-indices. Choosing a generating set $F_\nu$ of differential syzygies we express
 $$
\mathcal{A}=\la I_i\,;\nabla_j\,|\,F_\nu\ra.
 $$

This allows to solve the equivalence problem for submanifolds of functions with respect to $G$ as follows. Consider the above Lie-Tresse type
representation of $\mathcal{A}$. The collection of invariants $I_i,\nabla_j(I_i)$ (totally $r$ functions) allows to restore the generators,
while the relations $F_\nu$ constrain this collection. Any submanifold $N$ (for function $f$ given as the graph $\Sigma_f\simeq M$) canonically lifts
to the jet-space $J^\infty$: $N\ni a\mapsto [N]_a^\infty$. We thus map $\Psi:N\to\R^r$, $\Psi(a)=(I_i([N]_a^\infty),\nabla_j(I_i)([N]_a^\infty))$.
Due to differential syzygies the image is contained in some algebraic subset $Q\subset\R^r$. Two generic submanifolds $N_1,N_2$ are $G$-equivalent
iff $\Psi(N_1)=\Psi(N_2)$ as (un-parametrized) subsets.

\subsection{Conventions}
All differential invariants below are denoted by $I$ with a subscript. The subscript consists of a number and a letter.
The number reflects the order of an invariant, while the letter distinguishes invariants of the same order.
If no letter is given there is only one new (independent) invariant on the corresponding jet-space.

The symplectic Hamiltonian vector field in canonical coordinates on $V$ has the form
$X_H=\sum_i H_{y_i}\p_{x_i}-H_{x_i}\p_{y_i}$.
The Poisson bracket given by $[X_f,X_g]=X_{\{f,g\}}$ is equal to
 $$
\{f,g\}=\sum_{i=1}^n\left(\frac{\p f}{\p x_i}\frac{\p g}{\p y_i}-\frac{\p f}{\p y_i}\frac{\p g}{\p x_i}\right).
 $$
A basis of quadratic functions $\langle x_ix_j,x_iy_j,y_iy_j\rangle\ni f$ gives a basis of vector fields $X_f$ forming
$\sp(2n,\R)$. This may be extended to $\mathfrak{csp}(2n,\R)$ by adding the homothety $\zeta=\sum_i x_i\p_{x_i}+y_i\p_{y_i}$
that commutes with $\sp(2n,\R)$.

The contact Hamiltonian vector field in canonical coordinates on $W$ has the form
$X_H= H\p_z+ \sum_1^n\mathcal{D}^{(1)}_{x_i}(H)\p_{y_i}-H_{y_i}\mathcal{D}^{(1)}_{x_i}=
(H-\sum y_iH_{y_i})\p_z+\sum_1^n(H_{x_i}+y_iH_z)\p_{y_i}-H_{y_i}\p_{x_i}$.
The Lagrange bracket given by $[X_f,X_g]=X_{[f,g]}$ is equal to
 $$
[f,g]=\sum_{i=1}^n\left(\frac{\p f}{\p x_i}\frac{\p g}{\p y_i}-\frac{\p g}{\p x_i}\frac{\p f}{\p y_i}\right)
+\sum_{i=1}^n y_i\left(\frac{\p f}{\p z}\frac{\p g}{\p y_i}-\frac{\p g}{\p z}\frac{\p f}{\p y_i}\right)+\left(f\frac{\p g}{\p z}-g\frac{\p f}{\p z}\right).
 $$
A basis of quadratic functions $\langle x_ix_j,x_iy_j,y_iy_j\rangle\ni f$ gives a basis of vector fields $X_f$ forming
$\sp(2n,\R)$. This may be extended to $\mathfrak{csp}(2n,\R)$ by adding the homothety $X_f=\sum_i x_i\p_{x_i}+y_i\p_{y_i}+2z\p_z$ for $f=2z-\sum_ix_iy_i$
that commutes with $\sp(2n,\R)$.

\section{Functions on Symplectic Vector Spaces}

The group $G=\op{Sp}(2n,\R)$ acts almost transitively on $V=\R^{2n}$ (one open orbit that complements the origin),
and it is lifted to $J^0V=V\times\R(u)$ with $I_0=u$ being invariant. The prolonged action has
orbits of codimension 2 on $J^1V$ (one more invariant appears) and then the action becomes free on $J^2V$.

An invariant on $J^1$ is due to the invariant 1-form $du$ and the invariant (radial) vector field
$\zeta=\sum_i x_i\p_{x_i}+y_i\p_{y_i}$: their contraction yields
 $$
I_1=du(\zeta)=\sum_{i=1}^n x_i u_{x_i}+y_i u_{y_i}.
 $$

\subsection{The case of dimension $2n=2$}\label{F:n=1}

Here $V=\R^2(x,y)$.
To compute differential invariants of order $k$ we solve the equation
$\mathcal{L}_{X^{(k)}_i}I=0$, $I\in C^{\infty}(J^kV)$, for a basis
of the Lie algebra $\sp(2,\R)=\sl(2,\R)$: $X_1=x\p_y$, $X_2=x\p_x-y\p_y$, $X_3=y\p_x$.
For $k=2$, in addition to $I_0$ and $I_1$, we get
 \begin{align*}
I_{2a} &= x^2u_{xx}+2xyu_{xy}+y^2u_{yy}, \\
I_{2b} &= xu_yu_{xx}-yu_xu_{yy}+(yu_y-xu_x)u_{xy}, \\
I_{2c} &= u_x^2u_{yy}-2u_xu_yu_{xy}+u_y^2u_{xx}.
 \end{align*}
These invariants are functionally (hence algebraically) independent.

To determine the invariant derivations we solve its defining PDE. The invariant derivations of order $k=1$ are
linear combinations of
 $$
\nabla_1=x\mathcal{D}_x+y\mathcal{D}_y, \quad \nabla_2=u_x\mathcal{D}_y-u_y\mathcal{D}_x.
 $$

Let $\mathcal{A}$ denote the algebra of differential invariants, whose elements can be assumed polynomial in all jet-variables.
Since the obtained invariants are quasi-linear in their respective top jet-variables, and this property is
preserved by invariant derivations, the algebra $\mathcal{A}$ is generated by them.

To find a more compact description, note that $I_1=\nabla_1(I_0)$ and
 $$
I_{2a}=\nabla_{1}^{2}(I_{0})-\nabla_{1}(I_{0}), \quad
I_{2b}=-\nabla_{2}\nabla_{1}(I_{0}).
 $$
Thus only $I_0$ and $I_{2c}$ suffice to generate $\mathcal{A}$.

To describe the differential syzygies note that $\nabla_2(I_0)=0$, and the commutator relation is
 $$
[\nabla_{1},\nabla_{2}]=\frac{I_{2b}}{I_1}\nabla_1+\frac{I_{2a}-I_1}{I_1}\nabla_2.
 $$
Also, when applying $\nabla_1,\nabla_2$ to $I_{2a},I_{2b},I_{2c}$ and using the commutator relation
we get 5 different invariants of order 3, while there are only 4 independent 3-jet coordinates.
Thus computing the symbols of the invariants and eliminating those coordinates we obtain the remaining syzygy:
 $$
(\nabla_{2}(I_{2b})+\nabla_{1}(I_{2c}))I_{1}-(3I_{2a}-I_{1})I_{2c}+3I_{2b}^{2}=0.
 $$
To summarize, define
 \begin{align*}
\mathcal{R}_{1}&=\nabla_{2}(I_{0}), \\
\mathcal{R}_{2}&=I_1[\nabla_{1},\nabla_{2}]-I_{2b}\nabla_{1}-(I_{2a}-I_{1})\nabla_{2}, \\
\mathcal{R}_{3}&=I_1\nabla_{2}(I_{2b})+I_1\nabla_{1}(I_{2c})-(3I_{2a}-I_{1})I_{2c}+3I_{2b}^{2}.
 \end{align*}
Then, the algebra of differential invariants is given by generators and relations as follows:
\begin{framed}[.45\textwidth]
$$
\mathcal{A}=\la I_{0},I_{2c}\,;\,\nabla_{1},\nabla_{2}\ |\ \mathcal{R}_{1},\mathcal{R}_{2},\mathcal{R}_{3}\ra.$$
\end{framed}

\subsection{Another approach for $n=1$}
We act similar to \cite{KL1}.

Note that $\nabla_1$ corresponds to the radial vector field $\zeta$ and $\nabla_2=\omega^{-1}\hat d u$, where $\hat d$ is the horizontal differential
(in this case $\hat d=dx\otimes \mathcal{D}_x+dy\otimes \mathcal{D}_y$,
so $\hat d u=u_x\,dx+u_y\,dy$).
To find further invariants and derivations we consider the quadratic form
 $$
Q_2=d^2 u =u_{xx}dx^2+2u_{xy}dx\,dy+u_{yy}dy^2 \in \pi_2^*S^2T^*V.
 $$
Lowering the indices with respect to the symplectic form
(or partially contracting with $\omega^{-1}=\p_x\we\p_y$) we get the endomorphism
 $$
A=\omega^{-1}Q_2= u_{yy}\p_x\otimes dy-u_{xy}\p_y\otimes dy+u_{xy}\p_x\otimes dx-u_{xx} \p_y\otimes dx.
 $$
This can be lifted to the Cartan distribution on $J^\infty$ and thus applied to
horizontal fields:
 \begin{gather*}
A\nabla_1=(xu_{xy}+yu_{yy})\mathcal{D}_x-(xu_{xx}+yu_{xy})\mathcal{D}_y,\\
A\nabla_2=(u_xu_{yy}-u_yu_{xy})\mathcal{D}_x-(u_xu_{xy}-u_yu_{xx})\mathcal{D}_y.
 \end{gather*}
These are also invariant derivations and they can be expressed through the previous as follows:
 $$
A\nabla_1=-\frac{I_{2b}}{I_1}\nabla_1-\frac{I_{2a}}{I_1}\nabla_2,\quad
A\nabla_2=\frac{I_{2c}}{I_1}\nabla_1+\frac{I_{2b}}{I_1}\nabla_2.
 $$
Note also that $I_{2a}=Q_2(\nabla_1,\nabla_1)$, $I_{2b}=-Q_2(\nabla_1,\nabla_2)$, $I_{2c}=Q_2(\nabla_2,\nabla_2)$,
so that we can generate all the invariants uniformly.

\subsection{The general case}

In general dimension $2n$ we still have the invariant derivations
$\nabla_1$ corresponding to the radial field $\zeta$ and
$\nabla_2=\omega^{-1}Q_1$ for $Q_1=\hat d I_0$. Then the horizontal field of endomorphisms $A=\omega^{-1}Q_2$ for $Q_2=\hat d^2 I_0$
generates the rest: the invariant derivations $\nabla_{i+2}=A^i\nabla_2$
(alternatively $\nabla_{i+2}=A^i\nabla_1$) for $i=1,\dots,2n-2$
are independent (also with $\nabla_1,\nabla_2$)
on a Zariski open subset in the space of jets. This gives a
complete set of invariant derivations $\nabla_1,\dots,\nabla_{2n}$.

Taking into account $I_1=\nabla_1(I_0)$ the generating set of invariants
can be taken $I_0$ and $I_{ij}=Q_2(\nabla_i,\nabla_j)$.
By dimensional count and independence it is enough to
restrict to $i=1,2$ and $1\leq j\leq 2n$. We obtain:

 \begin{theorem}
The algebra of differential invariants of the $G$-action on $J^\infty(V)$ is
 $$
\mathcal{A}=\la I_{0},I_{1i},I_{2j}\,;\,\nabla_k\ |\ \mathcal{R}_l\ra
 $$
for some finite set of differential syzygies $\mathcal{R}_l$.
 \end{theorem}

This is a Lie-Tresse type of generation of $\mathcal{A}$.
Note also the following (non-finite) generation of this algebra.
The higher symmetric differentials $Q_k=d^ku \in \pi_k^*S^kT^*V$ can be
contracted with invariant derivations to get $k$-th order differential invariants
$Q_k(\nabla_{j_1},\dots,\nabla_{j_k})$. These clearly generate $\mathcal{A}$.

There is an algorithmic way of describing relations (syzygies)
between these invariants similar to \cite[Section 4]{KL1}.
We refer for explicit formulae of invariants to \cite{JJ} for $n=2$.

\section{Curves in Symplectic Vector Spaces}

Locally a curve in $\R^{2n}$ is given
as $\bu=\bu(t)$ for $t=x_1$ and $\bu=(x_2,\dots,x_n,y_1,\dots,y_n)$
in the canonical coordinates $(x_1,x_2,\dots,x_n,y_1,\dots,y_n)$ .
The corresponding jet-space $J^k(V,1)$ has coordinates $\bu_l$, $l\leq k$,
where $l$ stands for the $l$-tuple of $t$. For instance, $J^1(V,1)=\R^{4n-1}(t,\bu,\bu_1)$.
Note that $\dim J^k(V,1)=2n+k(2n-1)$.

\subsection{The case of dimension $2n=2$}\label{C:n=1}

Let us again start with the simplest example $V=\R^2(x,y)$.
The jet-space is $J^k(V,1)=\R^{k+2}(x,y,y_1,\dots,y_k)$.
Here $G=\op{Sp}(2,\R)$ has an open orbit in $J^1(V,1)$,
and there is one new differential invariant in every higher jet-order $k$.

Let us indicate in this simple case how to verify algebraicity of the action
(this easily generalizes to the other cases and will not be discussed further).
The 1-prolonged action of $g=\begin{pmatrix}a & b\\ c & d\end{pmatrix}\in G$ is
 $$
\Phi_{g}^{(1)}(x,y,y_{1})=\left(ax+by,cx+dy,\frac{dy_{1}+c}{by_{1}+a}\right).
 $$
Since the action on $J^0(V,1)=\R^2$ is transitive in the complement to 0,
we can choose point $p=(1,0)$ as a generic point. Its stabilizer group
is $G_p=\left\{\begin{pmatrix}1 & b\\ 0 & 1\end{pmatrix}\right\}\subset G$.
The action of this on the fiber $\pi_{1,0}^{-1}(p)$ is algebraic:
$y_1\mapsto\frac{y_{1}}{by_{1}+1}$. 

Thus the Lie-Tresse theorem \cite{KL} applies and the algebra of invariants
$\mathcal{A}$ can be taken to consist of rational functions in jet-variables,
which are polynomial in jets of order $\ge2$.

The first differential invariant is easily found from the Lie equation:
 $$
I_2=\frac{y_2}{(xy_1-y)^3}.
 $$
Similarly, solving the PDE for the coefficients of invariant derivation we find
 $$
\nabla=\frac1{xy_1-y}\mathcal{D}_x.
 $$
Now by differentiation we get new differential invariants
$I_3=\nabla I_2$, $I_4=\nabla^2 I_2$, etc. Since these are quasilinear
differential operators, they generate the entire algebra.
In other words, the algebra of differential invariants is free:
\begin{framed}[.2\textwidth]
 $$
\mathcal{A}=\la I_2\,;\,\nabla\ra.
 $$
\end{framed}

\subsection{The case of dimension $2n=4$}\label{C:n=2}

Let us use coordinates $(t,x,y,z)$ on $V=\R^4$ with the symplectic
form $\omega=dt\wedge dy+dx\wedge dz$.
Note that $\dim J^k(V,1)=3k+4$, and the jet-coordinates on $J^k$ are $(t,x,y,z,\dots,x_k,y_k,z_k)$. The action of $G=\op{Sp}(4,\R)$ on $J^k(V,1)$
has orbits of dimensions $4,7,9,10$ for $k=0,1,2,3$ respectively.
Thus the first differential invariant appears already in jet-order 2,
then two more appear in jet-order 3, and then there
appear $h_k=3$ new invariants in every jet-order $k\ge4$.

The infinitesimal and moving frame methods fail to produce enough invariants
here, so we apply more geometric considerations.

We exploit that $G$ preserves the symplectic form on $V$, but also the
fact that the action is linear, so the vector space structure of $V$ is preserved as well. In particular, the origin is preserved, so we can form a vector from the origin to any point $p=(t,x,y,z)\in J^{0}(V,1)$. Denote the corresponding vector by $$v_0=(t,x,y,z)\equiv t\p_t+x\p_x+y\p_y+z\p_z.$$

Consider the space of $1$-jets of unparametrized curves $J^1(V,1)$.
For a parameterization of the curve $c=(t,x(t),y(t),z(t))$ the
tangent vector at any point of this curve can be computed as $w_1=\mathcal{D}_t^{(1)}=\p_t+x_1\p_x+y_1\p_y+z_1\p_z$, which
is rescaled $v_1=\beta w_{1}$ upon a change of parametrization. To
make $v_1$ invariant we fix $\beta$ by the condition $\omega(v_{0},v_{1})=1$.
This normalization $\beta=1/(ty_1+xz_1-x_1z-y)$ gives
a canonical horizontal (that is tangent to the curve) vector field,
which can be interpreted as an invariant derivative
 $$
\nabla=\frac1{(ty_1+xz_1-x_1z-y)}\mathcal{D}_t.
 $$

The further approach is as follows. On every step there is a freedom
associated to a parameterization of a given curve. Fixing it
in a canonical way via evaluation with the symplectic form,
we obtain invariantly defined vectors and henceforth invariants.

On the first step, changing the parameterization $c=c(t)$
to another parameterization $c=c(\tau)$ results in a change of the
tangent vector by the chain rule:
 $$
\frac{dc}{dt}=\frac{d\tau}{dt}\frac{dc}{d\tau}.
 $$
This can be written as $w_1=k_1v_1$, for $d\tau/dt=k_1$.
The vector $w_1$, associated with a specific choice of parameterization,
is not canonical but convenient for computations.
The above normalization $k_1=1/\beta$ makes $v_1$ a canonical choice.

The change of parameterization on $2$-jets gives
 $$
\frac{d^2c}{dt^2}=\frac{d^2c}{d\tau^2}\left(\frac{d\tau}{dt}\right)^2
+\frac{dc}{d\tau}\frac{d^2\tau}{dt^2}.
 $$
Denote $v_2=d^2c/d\tau^2$, $w_2=d^2c/dt^2$ and $d^2\tau/dt^2=k_2$.
The equation becomes
 $$
w_{2}=v_{2}k_{1}^{2}+v_{1}k_{2}.
 $$
In the parameterization $c=c(t)$ the acceleration is $w_2=(0,x_2,y_2,z_2)$.
We solve for $v_2$ as
 $$
v_2=\frac{w_2-v_1k_2}{k_1^2}.
 $$
Then $k_2$ can be fixed by $\omega(v_0,v_2)=0$.
This uniquely determines $v_2$, which can now be used to find the first differential invariant. In fact, $I_2=\omega(v_1,v_2)$ is a differential invariant of order $2$. In coordinates
 $$
I_2=\omega(v_1,v_2)=\frac{x_1z_2-z_1x_2+y_2}{(ty_1+xz_1-zx_1-y)^3}.
 $$

There are $2$ independent third order invariants by our dimension count.
The first can be obtained as $\nabla(I_2)$, to find the second we exploit
the above normalization method on 3-jets. The change of parameterization is
 $$
\frac{d^3c}{dt^3}=\frac{d^3c}{d\tau^3}\left(\frac{d\tau}{dt}\right)^3
+3\frac{d^2c}{d\tau^2}\frac{d\tau}{dt}\frac{d^2\tau}{dt^2}
+\frac{dc}{d\tau}\frac{d^3\tau}{dt^3}.
 $$
Again, rewrite it in simpler notations as
 $$
w_3 = v_3k_1^3+3k_1k_2v_2+k_3v_1.
 $$
Here $w_3=(0,x_3,y_3,z_3)$ and the unknown $k_3$ can be fixed by
the condition $\omega(v_0,v_3)=0$, where
 $$
v_3=\frac{w_3-3k_1k_2v_2-k_3v_1}{k_1^3}.
 $$
This uniquely determines $v_3$, which allows the computation of two new
differential invariants:
 $$
I_{3a}= \omega(v_1,v_3), \quad I_{3b}= \omega(v_2,v_3).
 $$
The invariants $I_{3a}$ and $I_{3b}$ are independent, but $I_{3a}$ can be expressed through $\nabla(I_{2})$, so it is not required in what follows.

Finally we explore the forth order chain rule
 $$
\frac{d^4c}{dt^4}=\frac{d^4c}{d\tau^4}\left(\frac{d\tau}{dt}\right)^4
+6\frac{d^3c}{d\tau^3}\left(\frac{d\tau}{dt}\right)^2\frac{d^2\tau}{dt^2}
+\frac{d^2c}{d\tau^2}\left(4\frac{d\tau}{dt}\frac{d^3\tau}{dt^3}
+3\left(\frac{d^2\tau}{dt^2}\right)^2\right)
+\frac{dc}{d\tau}\frac{d^4\tau}{dt^4}
 $$
that can be written as
 $$
w_4= v_{4}k_{1}^{4}+6v_{3}k_{1}^{2}k_{2}+v_{2}\left(4k_{1}k_{3}+3k_{2}^{2}\right)+v_{1}k_{4}
 $$
with $w_4=(0,x_4,y_4,z_4)$. Find $k_{4}$ by $\omega(v_{0},v_{4})=0$.
This uniquely determines $v_4$, then the invariants of order $4$
are found by the formulae
 $$
I_{4a}=\omega(v_{1},v_{4}), \quad
I_{4b}=\omega(v_{2},v_{4}), \quad
I_{4c}=\omega(v_{3},v_{4}).
 $$
These are independent, but $I_{4a}$ and $I_{4b}$ can be expressed by the invariants of order $3$ and the invariant derivation, so they
will not be required in what follows.

This gives the necessary invariants to generate the entire algebra of differential invariants. To summarize, if we denote $I_3=I_{3b}$ and
$I_4=I_{4c}$, then the algebra of differential invariants is freely generated as follows
\begin{framed}[.3\textwidth]
 $$
\mathcal{A}=\la I_2,I_3,I_4\,;\,\nabla\ra.
 $$
\end{framed}
The explicit coordinate formulae of invariants are shown in the Appendix.

\subsection{The general case}

In dimension $\dim V=2n$ the following dimensional analysis
readily follows from the normalization procedure developed above.

\smallskip

 \begin{center}
\begin{tabular}{|l|l|l|l|} \hline
Jet-level $k$ & $\dim J^k(V,1)$ & $G$-orbit dimension & \# new invariants $h_k$ \\ \hline
$0$ & $2n$ & $2n$ & $0$ \\
$1$ & $4n-1$ & $2n+(2n-1)=4n-1$ & $0$ \\
$2$ & $6n-2$ & $(4n-1)+(2n-2)=6n-3$ & $1$ \\
$3$ & $8n-3$ & $(6n-3)+(2n-3)=8n-6$ & $2$ \\
$4$ & $10n-4$ & $(8n-6)+(2n-4)=10n-10$ & $3$ \\
\dots & \dots & \dots & \dots \\
$k$ & $2n+k(2n-1)$ & $2(k+1)n-\binom{k+1}2$ & $k-1$ \\
\dots & \dots & \dots & \dots \\
$2n-1$ & $(2n-1)^2+2n$ & $\binom{2n+1}2$ & $2n-2$ \\
$2n$ & $4n^2$ & $\dashrightarrow$ stabilized & $2n-1$\\ \hline
\end{tabular}
 \end{center}

\smallskip

In particular, the number of pure order $k$ differential invariants
is $h_k=k-1$ for $1\leq k\leq 2n$ and $h_k=2n-1$ for $k>2n$.

If the canonical coordinates in $\R^{2n}$ are $(t,\bx,y,\bz)$,
where $\bx$ and $\bz$ and $(n-1)$-dimensional vectors, then
the invariant derivation is equal to
 $$
\nabla=\frac1{(ty_1-y+\bx \bz_1-\bx_1 \bz)}\mathcal{D}_t.
 $$
We also obtain the first differential invariant of order 2
 $$
I_2=\frac{(\bx_1 \bz_2-\bx_2 \bz_1+y_2)}{(ty_1-y+\bx \bz_1-\bx_1 \bz)^3}.
 $$
Then we derive the differential invariant $\nabla(I_2)$
and add to it another differential invariant $I_3$ of order 3.
Then we derive the differential invariants $\nabla^2(I_2),\nabla(I_3)$
and add another differential invariant $I_4$ of order 4.
We continue obtaining new invariants by using the higher order chain rule
and normalization via the symplectic form up to order $2n$.

In summary, we obtain $2n-1$ independent
differential invariants $I_2,\dots,I_{2n}$
of orders $2,\dots,2n$ respectively.

 \begin{theorem}
The algebra of differential invariants of the $G$-action on $J^\infty(V,1)$
is freely generated as follows:
 $$
\mathcal{A}=\la I_2,\dots,I_{2n}\,;\, \nabla\ra.
 $$
 \end{theorem}

\section{Hypersurfaces in Symplectic Vector Spaces}

Since hypersurfaces in $\R^2$ are curves,
the first new case come in dimension 4. We consider
this first and then discuss the general case.

\subsection{The case of dimension $2n=4$.}

Let $V=\mathbb{R}^{4}$, denote its canonical coordinates by $(x,y,z,u)$ with $\omega=dx\wedge dz+dy\wedge du$.
Hypersurfaces can be locally identified as graphs $u=u(x,y,z)$ and this gives parametrization of an open chart in $J^k(V,3)$. We
use the usual jet-coordinates $u_x$, $u_{xx}$, $u_{yz}$, etc.

As is the cases above, straightforward computations become harder. Maple is not able to compute all required invariants and derivations,
so we again rely on a more geometric approach.
Before going through the method we investigate the count of invariants.

The group $G=\op{Sp}(4,\R)$ acts with an open orbit on $J^0(V,3)$. On the space of 1-jets the dimension of the orbit is $7=\dim J^1(V,3)$,
hence there are no invariants. The orbit stabilization is reached on $J^2(V,3)$, where the action is free. The rank of the action is $10$ and
$\dim J^2(V,3)=13$, so there are $h_2=3$ independent differential invariants. For $k>2$, the number of new differential invariants is
$h_k=\binom{k+2}2$. In particular, $h_3=10$.

The number of independent invariant derivations is $3$, so these and 3 invariants of order 2 generate a total number of $9$ invariants of order 3.
In addition, commutators of invariant derivations $[\nabla_i,\nabla_j]=I_{ij}^k\nabla_k$ give up to $9$ more differential invariants of order~$3$.
We will confirm that the totality of these $18$ contain $10$ independent invariants of order~3, and hence suffice to generate
also the differential invariants of higher order.

The 0-jet $p=(x,y,z,u)\in J^0(V,3)$ can be identified with the vector from the origin to this point, which we denote by
$$v_0=(x,y,z,u)\equiv x\p_x+y\p_y+z\p_z+u\p_u.$$
The 1-jet of a hypersurface $\Sigma=\{u=u(x,y,z)\}$ can be identified with its tangent space
 $$
T_p\Sigma= \la\p_x+u_x\p_u,\p_y+u_y\p_u,\p_z+u_z\p_u\ra= \la\mathcal{D}_x^{(1)},\mathcal{D}_y^{(1)},\mathcal{D}_z^{(1)}\ra.
 $$
The orthogonal complement to $T_p\Sigma$ with respect to $\omega$ is generated by
 $$
w_1=\p_y-u_z\p_x+u_x\p_z+u_y\p_u,
 $$
that is $T_p\Sigma^{\perp\omega}=\la w_1\ra$. The vector $w_1$ is determined up to scale, which we fix via the symplectic form so:
$v_1=k_1w_1$ must satisfy $\omega(v_0,v_1)=1$. This normalization gives $k_1=1/(xu_x+yu_y+zu_z-u)$, so the canonical vector $v_1$ is equal to
 $$
v_1 = \frac1{xu_x+yu_y+zu_z-u}(\partial_y-u_z\partial_x+u_x\partial_z+u_y\partial_u).
 $$
This vector field is tangent to the hypersurface, so it is horizontal and can be rewritten in terms of the total derivative.
This yields the first invariant derivation:
 $$
\nabla_1=\frac{\mathcal{D}_y-u_z\mathcal{D}_x+u_x\mathcal{D}_z}{xu_x+yu_y+zu_z-u}.
 $$

Let $q=-u+u(x,y,z)$ be a defining function of the hypersurface $\Sigma=\{q=0\}$. We have $T_p\Sigma=\op{Ker}{dq}$.
A change of the defining function $q'=fq$ of $\Sigma$, with $f\in C^\infty(V)$ such that $f|_{\Sigma}\neq 0$, has the following effect on
the differential: $dq'=q\,df+f\,dq$. Therefore at $p\in\Sigma$ we have $d_pq'=f(p)d_pq$ and so $T_p\Sigma=\op{Ker}{dq'}$.

Next we compute the second symmetric differential $d^2q$ of the defining function for $\Sigma$. A change of the defining function $q'=fq$
has the following effect on the second differential:
 $$
d^2q'=d(d(fq))=d(q\,df+f\,dq)=q\,d^2f+2\,df\,dq+f\,d^2q.
 $$
At the points $p\in\Sigma$ this simplifies to
 $$
d_p^2q'=2\,d_pf\,d_pq+f(p)\,d_p^2q.
 $$
Restricting to the tangent space of $\Sigma$ gives
 $$
d^2q'\big|_{T_p\Sigma}=f(p)d^2q\big|_{T_p\Sigma}.
 $$

Thus the defining differential $dq$ and the quadratic form $d^2q$ are defined up to the same scale.
We fix it again via the symplectic form: $d_pq'=k_2d_pq$ must satisfy $d_pq'(v_0)=1$, i.e.\ $k_{2}=1/dq(v_{0})$ for generic $1$-jets.
This normalization gives the quadratic form $d^2q'|_{T_p\Sigma}=k_2d^2q\big|_{T_p\Sigma}=d^2q\big|_{T_p\Sigma}/dq(v_0)$.
In coordinates, with $q=-u+u(x,y,z)$, we get the expression
 $$
Q=d^{2}q'\big|_{T_p\Sigma}=\frac{u_{xx}dx^2+2u_{xy}dxdy+2u_{xz}dxdz+u_{yy}dy^2+2u_{yz}dydz+u_{zz}dz^2}{xu_x+yu_y+zu_z-u}.
 $$
The first invariant is then computed by
 $$
I_{2a}=Q(v_1,v_1)=\frac{u_x^2u_{zz}-2u_xu_zu_{xz}+u_z^2u_{xx}+2u_xu_{yz}-2u_zu_{xy}+u_{yy}}{(xu_x+yu_y+zu_z-u)^3}.
 $$

Let us summarize the geometric data encoding the 2-jet that we obtained and which are supported on the 3-dimensional tangent space $T_p\Sigma$:
the invariant vector $v_1$, the symmetric 2-form $Q$ of general rank, the skew 2-form $\omega|_{T_p\Sigma}$ of rank $2$ ($v_1$ spans its kernel),
and $1$-form $\alpha=\omega(v_0,\cdot)$. These data give a canonical splitting of the tangent space $T_p\Sigma=\la v_1\ra\oplus\Pi$,
where $\Pi=\op{Ker}(\alpha)$. Indeed, $v_1\notin\op{Ker}(\alpha)$ because $\omega(v_0,v_1)=1$ by the normalization.
Using this data we can construct 2 more invariant derivations.

Choose a nonzero $w_3\in\Pi$, $Q(v_1,w_3)=0$. Then choose $w_2\in\Pi$, $Q(w_2,w_3)=0$.
For generic 2-jet, the vectors $w_2,w_3$ are defined up to scale that we fix so: $v_2\in\la w_2\ra$, $v_3\in\la w_3\ra$
must satisfy $Q(v_1,v_2)=1$, $\omega(v_2,v_3)=1$.

Since $v_2,v_3\in T_p\Sigma$ are horizontal, they generate two invariant derivations $\nabla_2,\nabla_3$.
Additionally we get 2 differential invariants:
 $$
I_{2b}=Q(v_2,v_2),\quad I_{2c}=Q(v_3,v_3).
 $$
A calculation of the rank of the corresponding Jacobi matrix shows that these are independent, and moreover that the data generate all differential
invariants of order 3. Then by independence of $\nabla_1,\nabla_2,\nabla_3$ all higher order invariants can be derived, so
for a finite set of differential syzygies $\mathcal{R}_l$ we get:
\begin{framed}[.5\textwidth]
 $$
\mathcal{A}=\la I_{2a},I_{2b},I_{2c}\,;\,\nabla_1,\nabla_2,\nabla_3 \ |\ \mathcal{R}_l\ra
 $$
\end{framed}

The coordinate formulae can be found in \cite{JJ} (note that renumeration $v_2\leftrightarrow v_3$ and
a different normalization is taken here).

\subsection{The general case}

Now we consider jets of hypersurfaces $\Sigma\subset V=\R^{2n}$ for general $n$ and compute their differential invariants
with respect to $G=\op{Sp}(2n,\R)$.

By the Lie-Tresse theorem \cite{KL} the algebra $\mathcal{A}$ can be assumed to consist of rational functions on $J^\infty(V,2n-1)$,
which are polynomial in jet-variables of order $\ge 2$.

The dimensional count easily generalizes to give $h_0=h_1=0$, $h_2=2n-1$ and $h_k=\binom{2n-2+k}{k}$ for $k>2$.
There will be $2n-1$ independent invariant derivations $\nabla_j$, and as before these together with second order invariants $I_{2s}$
($1\leq s\leq 2n-1$) and the structure coefficients $I_{ij}^k$ of the horizontal frame $\nabla_j$ will suffice to generate all invariants.

We again have the position vector $v_0$, the tangent vector $v_1$ normalized by $\omega(v_0,v_1)=1$, and
the quadratic form $Q$ on $T_p\Sigma$. From this data in a Zariski open set of $J^2(V,2n-1)$ of generic 2-jets we get a canonical basis
$v_1,\dots,v_{2n-1}$ by normalizing in turn via $\omega$ and $Q$ as follows (we repeat steps 0 and 1 that are already performed).

Step 0: $T_p\Sigma=\la v_1,\dots,v_{2n-1}\ra$.

Step 1: Choose $v_1$ by $\la v_1\ra\perp_\omega\la v_1,\dots,v_{2n-1}\ra$, $\la v_2,\dots,v_{2n-1}\ra\perp_\omega\la v_0\ra$.
Normalize $\omega(v_0,v_1)=1$.

Step 2: Choose $v_2$ by $\la v_3,\dots,v_{2n-1}\ra\perp_Q\la v_1\ra$, $\la v_2\ra\perp_Q\la v_3,\dots,v_{2n-1}\ra$.
Normalize $Q(v_1,v_2)=1$.

Step 3: Choose $v_3$ by $\la v_3\ra\perp_\omega\la v_3,\dots,v_{2n-1}\ra$, $\la v_4,\dots,v_{2n-1}\ra\perp_\omega\la v_2\ra$.
Normalize $\omega(v_2,v_3)=1$.

Step 4: Choose $v_4$ by $\la v_5,\dots,v_{2n-1}\ra\perp_Q\la v_3\ra$, $\la v_4\ra\perp_Q\la v_5,\dots,v_{2n-1}\ra$.
Normalize $Q(v_3,v_4)=1$.

Inductively we get the interchangeable steps as follows.

Step $(2r-1)$: Choose $v_{2r-1}$ by $\la v_{2r-1}\ra\perp_\omega\la v_{2r-1},\dots,v_{2n-1}\ra$, $\la v_{2r},\dots,v_{2n-1}\ra\perp_\omega\la v_{2r-2}\ra$.
Normalize $\omega(v_{2r-2},v_{2r-1})=1$.

Step $2r$: Choose $v_{2r}$ by $\la v_{2r+1},\dots,v_{2n-1}\ra\perp_Q\la v_{2r-1}\ra$, $\la v_{2r}\ra\perp_Q\la v_{2r+1},\dots,v_{2n-1}\ra$.
Normalize $Q(v_{2r-1},v_{2r})=1$.

The procedure stops at step $(2n-1)$. The frame $v_i$ is canonical:
 $$
\omega^{-1}=v_0\wedge v_1+v_2\wedge v_3+\dots+v_{2n-2}\wedge v_{2n-1}.
 $$
The only non-constant entries of the Gram matrix of $Q$ in the basis $v_i$
are diagonal $Q(v_i,v_i)=I_{2,i}$ for $1\leq i<2n$.
The Gram matrix consists of $(n-1)$ diagonal blocks of size $2\times2$
and 1 diagonal block of size $1\times1$ as follows:
 \begin{center}
\begin{tabular}{c|cccccccc}
$Q$ & $v_1$ & $v_2$ & $v_3$ & $v_4$ & \dots & $v_{2n-3}$ & $v_{2n-2}$ & $v_{2n-1}$\\ \hline
$v_1$ & $I_{2,1}$ & 1 & 0 & 0 & \dots & 0 & 0 & 0 \\
$v_2$ & 1 & $I_{2,2}$ & 0 & 0 & \dots & 0 & 0 & 0 \\
$v_3$ & 0 & 0 & $I_{2,3}$ & 1 & \dots & 0 & 0 & 0 \\
$v_4$ & 0 & 0 & 1 & $I_{2,4}$ & \dots & 0 & 0 & 0 \\
$\vdots$ & $\vdots$ & $\vdots$ & $\vdots$ & $\vdots$ & $\ddots$ & $\vdots$ & $\vdots$ & $\vdots$ \\
$v_{2n-3}$ & 0 & 0 & 0 & 0 & \dots & $I_{2,2n-3}$ & 1 & 0 \\
$v_{2n-2}$ & 0 & 0 & 0 & 0 & \dots & 1 & $I_{2,2n-2}$ & 0 \\
$v_{2n-1}$ & 0 & 0 & 0 & 0 & \dots & 0 & 0 & $I_{2,2n-1}$
\end{tabular}
 \end{center}

The horizontal vector fields $v_j$ correspond to invariant derivations
$\nabla_j$, $1\leq j\leq 2n-1$. To summarize we obtain the following statement.

 \begin{theorem}
For the $G$-action on $J^\infty(V,2n-1)$ the algebra $\mathcal{A}$ 
is generated by the differential invariants
$I_{2,i}$ and the invariant derivations $\nabla_j$, where $1\leq i,j\leq 2n-1$.
 \end{theorem}

\section{General submanifolds in a Symplectic Vector Space}

The case of submanifolds of dimension and codimension greater than 1 is
more complicated, no straightforward computations work
for $G=\op{Sp}(2n,\R)$ action on $J^\infty(V,m)$, $V=\R^{2n}$.
Yet the geometric methods applied above do generalize, and to
illustrate this we consider the simplest case $n=m=2$ and then remark
on the general case.

\subsection{Surfaces in a four-dimensional Symplectic Space}

The action has an open orbit in $J^1(V,2)$, but becomes free
on the level of 2-jets. Since $\dim J^2(V,2)=14$ we get
$h_2=4$ differential invariants of order 2 and then
at every higher order $k>2$ there will be $h_k=2(k+1)$ new invariants.

There will be 2 independent invariant derivations. Applying those
to 4 differential invariants of the second order gives a total of $8$
invariants of order 3. A direct computation shows that these are
functionally (hence algebraic) independent.
Since $h_3=8$ this is enough to generate all differential invariants.

In this case 
the algebra $\mathcal{A}$ of differential invariants can be chosen to consist of rational functions that are polynomial in jets-variables of order $>2$.

Having done the counting we can proceed with the geometric approach.
Choose canonical coordinates $(t,s,x,y)$ on $V=(\R^4,\omega)$.
Locally surfaces in $V$ are given as
$\Sigma=\{x=x(s,t),y=y(s,t)\}$. Here $s,t$ will be treated as independent
and $x,y$ as dependent variables, whence the coordinates on $J^\infty(V,2)$.

The 0-jet $p=(t,s,x,y)\in J^0$ can be identified with the vector
to that point from the origin
$v_0=t\p_t+s\p_s+x\p_x+y\p_y$.The 1-jet can be identified with
the tangent space
 $$
T_p\Sigma=\la\mathcal{D}_t^{(1)},\mathcal{D}_s^{(1)}\ra=
\la\p_t+x_t\p_x+y_t\p_y,\p_s+x_s\p_x+y_s\p_y\ra.
 $$
Equivalently, if the surface is described by $\Sigma=\{f=0,g=0\}$ with $f=x-x(t,s)$ and $g=y-y(t,s)$, then $T_p\Sigma=\text{Ann}(d_pf,d_pg)$,
where $d_pf=dx-x_tdt-x_sds$ and $d_pg=dy-y_tdt-y_sds$.

The restriction of $\omega$ to $T_p\Sigma$ has rank $2$ on generic $1$-jets, so $T_p\Sigma$ is a symplectic subspace of dimension $2$ and
$T_pV=T_p\Sigma\oplus T_p\Sigma^{\perp\omega}$.

Denote by $\pi_1:T_pV\to T_p\Sigma$ and $\pi_2:T_pV\to T_p\Sigma^{\perp\omega}$
the natural projections with respect to this decomposition.
Further for $v\in T_pV$ denote $v=v^{\parallel}+v^{\perp}$, where $v^{\parallel}=\pi_1(v)\in T_p\Sigma$ and $v^{\perp}=\pi_2(v)\in T_p\Sigma
^{\perp\omega}$.

Thus, 1-jet $[\Sigma]^1_p$ is entirely encoded by
$(T_p\Sigma,\omega|_{T_p\Sigma},v_0^{\parallel})$ and
$(T_p\Sigma^{\perp\omega},\omega|_{T_p\Sigma^{\perp\omega}},v_0^{\perp})$.
Note also that $\op{Ann}(T_p\Sigma)$ is identified with
$T_p\Sigma^{\perp\omega}$ by the symplectic form $\omega$.

Moving on to $2$-jets there is more structure on the tangent space.
The defining functions $f$, $g$ can be changed to $F=\alpha f+\beta g$,
$G=\gamma f+\delta g$, where $\alpha,\beta,\gamma,\delta$ are arbitrary functions that satisfy $\alpha\delta-\beta\gamma\neq 0$ along $\Sigma$.
Then $\Sigma=\{F=0,G=0\}$ and the tangent space can be described
as the annihilator of the differentials of the new defining functions
at $p\in\Sigma$:
 \begin{align*}
d_{p}F &= \alpha(p)d_{p}f+\beta(p)d_{p}g, \\
d_{p}G &= \gamma(p)d_{p}f+\delta(p)d_{p}g.
 \end{align*}
Next, compute the second symmetric differential of $f,g$ and
restrict to $T_p\Sigma$. Doing the same for $F,G$ results in
 \begin{align*}
d^2_pF &= \alpha(p)d^2_pf+\beta(p)d^2_pg, \\
d^2_pG &= \gamma(p)d^2_pf+\delta(p)d^2_pg.
 \end{align*}
This gives a $2$-dimensional space
$\mathcal{Q}=\la d_p^2f|_{T_p\Sigma},d_p^2g|_{T_p\Sigma}\ra=
\la d_p^2F|_{T_p\Sigma},d_p^2G|_{T_p\Sigma}\ra$ of quadratic forms,
and the above formulae show that there is a natural isomorphism
between $\op{Ann}(T_p\Sigma)\subset T_p^*V$ and $\mathcal{Q}$.
Our goal is to find a canonical basis $Q_1,Q_2$ in this space.

Let $Q_1\in\mathcal{Q}$ be given by the condition $Q_1(v_0^{\parallel},v_0^{\parallel})=0$. This ensures that $Q_1$
has a Lorentzian signature or is degenerate, and for a generic 2-jet
we get that $Q_1$ is non-degenerate. The vector $v_0^{\parallel}$ becomes
null-like vector for $Q_1$ that is yet defined up to scale.
A Lorentzian metric on the plane has two independent null-like vectors
and this gives a way to fix $Q_1$ and a vector $w^{\parallel}\in T_p\Sigma$
complementary to $v_0^{\parallel}$ as follows:
 $$
\omega(v_0^{\parallel},w^{\parallel})=1,\
Q_1(w^{\parallel},w^{\parallel})=0,\ Q_1(v_0^{\parallel},w^{\parallel})=1.
 $$
Note that this does not involve square roots, but only linear algebra.
Indeed, the first condition fixes the second null-like vector up to change
$w^{\parallel}\mapsto w^{\parallel}+kv_{0}^{\parallel}$.
The second condition fixes $k$ and the last normalizes $Q_1$.

The quadratic form $Q_1$ corresponds to a $1$-form $\sigma_1\in\text{Ann}(T_p\Sigma)$ such that the symmetric differential
of an extension of $\sigma_1$ to a section of $\text{Ann}(T\Sigma)$,
restricted to $T_p\Sigma$ equals $Q_1=d^\text{sym}_p\sigma_1$.
Then fix $w^\perp\in T_p\Sigma^{\perp\omega}$ uniquely by the conditions
$\sigma_1(w^{\perp})=0$, $\omega(v_0^{\perp},w^{\perp})=1$
(for a generic 2-jet $\sigma_1(v_0^{\perp})\neq0$).

Then define $\sigma_2\in\text{Ann}(T_p\Sigma)$ by the conditions $\sigma_2(v_0^{\perp})=0$, $\sigma_2(w^{\perp})=1$. This gives a
unique $1$-form independent of $\sigma_1$. It in turn
corresponds to a quadratic form $Q_{2}=d_p^\text{sym}\sigma_2$.

The remaining evaluations yield differential invariants
 $$
I_{2a}=\sigma_{1}(v_{0}^{\perp}),\ \
I_{2b}= Q_{2}(v_{0}^{\parallel},v_{0}^{\parallel}), \ \
I_{2c}= Q_{2}(v_{0}^{\parallel},w^{\parallel}), \ \
I_{2d}= Q_{2}(w^{\parallel},w^{\parallel}).
 $$
The vectors $v_0^{\parallel}$ and $w^{\parallel}$ are tangent vectors to $\Sigma$ (horizontal) so they correspond to the invariant derivations
$\nabla_1,\nabla_2$ and we conclude:
 \begin{theorem}
For the $G$-action on $J^\infty(V,2)$ the algebra $\mathcal{A}$
is generated by the differential invariants
$I_{2a},I_{2b},I_{2c},I_{2d}$ and invariant derivations $\nabla_1,\nabla_2$.
 \end{theorem}

Explicit form of these generators in jet-coordinates can be found in
\cite{JJ}.

\subsection{A remark about the general case}

In general it is easy to check that $G=\op{Sp}(2n,\R)$ acts
with one open orbit in $J^1(V,m)$, $V=\R^{2n}$, so there are no
first order invariants. However there are always second order invariants.
Their number is at least $\dim J^2(V,m)-n(2n+1)$, but this can be
non-positive for $m\ll n$.

Thus combining the ideas on differentials and quadratic forms with
$\omega$-orthogonal complements one can get some of the invariants.
If they are not sufficient, third and higher symmetric powers $d^rf$
of the defining functions $f$ should be explored.

From the investigated cases we cannot observe a pattern and hence
cannot universally describe all differential invariants of
$G=\op{Sp}(2n,\R)$ action on $J^\infty(V,m)$, $V=\R^{2n}$.

\section{Note on extension of the group}\label{S7}

One can also consider invariants of functions and submanifolds
in symplectic $V=\R^{2n}$ with respect to conformal symplectic
group $\op{CSp}(2n,\R)=\op{Sp}(2n,\R)\times\R_+$, the affine
symplectic group $\op{ASp}(2n,\R)=\op{Sp}(2n,\R)\ltimes\R^{2n}$
and affine conformal symplectic group
$\op{ACSp}(2n,\R)=\op{CSp}(2n,\R)\ltimes\R^{2n}$.
Denote a group in this list by $H$.

Since our $G$ is a subgroup of $H$, the algebras of differential
invariants $\mathcal{A}_H$ for each of the cases are subalgebras in
the algebra $\mathcal{A}_G$ that we previously computed
(enhanced notations should be self-obvious).
One imposes the homogeneity assumption or translation-invariance
or both on a general combination of invariants.

Let us discuss how to do this in all three cases.
For brevity of exposition we restrict to the case $n=1$
(functions and curves on symplectic plane), the general case is similar.

\subsection{Conformal symplectic group action: functions}\label{s71}

Consider functions on the conformal symplectic plane, 
$H=\op{CSp}(2n,\R)$. For $n=1$ observe $H\simeq\op{GL}(2,\R)$.
We recall the invariants from Section \ref{F:n=1} and note that
all of them are homogeneous with respect to scaling $\xi=x\p_x+y\p_y$,
corresponding to the center of $\mathfrak{h}=\gl(2,\R)$.
Restricting to invariants and derivations of weight 0 we
obtain the algebra of differential $\mathfrak{h}$-invariants.

The invariants $I_0,I_1$ have weight 0, and the invariants $I_{2a},I_{2b},I_{2c}$
have weights $0,-2,-4$ respectively. Therefore for the new algebra
$\mathcal{A}_H$ there are two independent invariants of order $\leq1$
and two additional invariants of order 2, namely $I_0$, $I_1$, $I_{2a}$
and $I'_{2b}=I_{2b}^{-2}I_{2c}$ in the notations of \S\ref{F:n=1}.

The invariant derivations are $\nabla_1,\nabla_2$ of weights $0,-2$
respectively. Therefore we obtain two invariant derivations with
respect to $\mathfrak{h}$: $\nabla_1$ and $\nabla'_2=I_{2b}^{-1}\nabla_2$.

Now a straightforward verification shows that $\nabla_1(I_{2a}),\nabla_1(I'_{2b}),
\nabla'_2(I_{2a}),\nabla'_2(I'_{2b})$ are independent in 3-jets, which implies that
the algebra $\mathcal{A}_H$ of differential invariants
is generated by $I_0,I'_{2b}$ and $\nabla_1,\nabla'_2$.
Note that $I_1=\nabla_1(I_0)$ and $I_{2a}=\nabla_1(I_1)-I_1$. 

To complete the picture, here are the differential syzygies: $\nabla'_2(I_0)=0$, $\nabla'_2(I_1)=-1$ and
 $$
[\nabla_1,\nabla'_2]=\frac1{I_1}\nabla_1+\Bigl(\frac{I_{2a}}{I_1}+\nabla'_2(I_{2a})\Bigr)\nabla'_2.
 $$
Denote these by $\mathcal{R}_1,\mathcal{R}_2,\mathcal{R}_3$.
There is also a forth order differential syzygy $\mathcal{R}_4$:
\begin{multline*}
\nabla_{2}'(I_{3a})+\frac{1}{2I'_{2b}}\nabla_{2}'(I_{3b})
-\frac{1}{2I'_{2b}}\nabla_{1}(I_{3c})-\frac{1}{2I_{1}^{2}I'_{2b}}
\Bigl[(I'_{2b}I_{3b}-3I_{3b}I_{3c}-I_{3c})I_{1}^{2}\\
+\left(\bigl((3I_{3b}+4)I_{2a}-5I_{3a}\bigr)I'_{2b}-4I_{2a}I_{3c}-4I_{3b}-4\right)I_{1}
+6I_{2a}^{2}I'_{2b}-6I_{2a}\Bigr]=0,
\end{multline*}
where $I_{3a}=\nabla_{1}(I_{2a}), I_{3b}=\nabla'_{2}(I_{2a}), I_{3c}=\nabla_{1}(I'_{2b})$ and $I_{3d}=\nabla'_{2}(I'_{2b})$. With this we obtain a complete description of the algebra
of differential $H$-invariants:
\begin{framed}[.5\textwidth]
 $$
\mathcal{A}_H=\la I_0,I'_{2b}\,;\,\nabla_1,\nabla'_2\,|\,
\mathcal{R}_1,\mathcal{R}_2,\mathcal{R}_3,\mathcal{R}_4\ra.
 $$
\end{framed}

\subsection{Conformal symplectic group action: curves}

Now we discuss differential invariants of curves with respect
to the same $H$ as in \S\ref{s71}.
Consider the invariants from Section \ref{C:n=1} and note that
all of them are homogeneous with respect to scaling $\xi=x\p_x+y\p_y$,
corresponding to the center of $\mathfrak{h}$.
Again we have to restrict to invariants and derivations of weight 0
to describe the algebra $\mathcal{A}_H$.

The invariant $I_2$ has weight $-4$ and the derivation $\nabla$
weight $-2$. Thus the derived invariants $I_{k+2}=\nabla^k(I_2)$ have weights
$-2(k+2)$ for $k\ge0$. In particular, $I'_3=I_3^2/I_2^3$ has weight 0
and similar for $\nabla'=I_2I_3^{-1}\nabla$ in the notations of \S\ref{C:n=1}.
Therefore these freely generate the algebra of differential $H$-invariants:
\begin{framed}[.2\textwidth]
  $$
\mathcal{A}_H=\la I'_3\,;\,\nabla'\ra.
 $$
\end{framed}

\subsection{Affine symplectic group action: functions}\label{asga1}

Consider differential invariants of functions on the
affine symplectic plane, $H=\op{ASp}(2n,\R)$.
For $n=1$ observe $H=\op{SAff}(2,\R)$.
We recall the generating invariants from Section \ref{F:n=1}, and note
that they indeed depend explicitly on $x,y$ except for $I_0$ and $I_{2c}$.

To single out invariants in $\mathcal{A}_G$ that are $x,y$-independent
eliminate $x,y$ from the system $\{I_1=c_1,I_{2a}=c_2,I_{2b}=c_3\}$
to get a translation-invariant polynomial on $J^2$ that depends parametrically
on $c_1,c_2,c_3$. Taking the coefficients of this expression with respect to
those parameters, we obtain the invariants $I_{2c}$ and
$I'_2=u_{xx}u_{yy}-u_{xy}^2=\op{Hess}(u)$. Then substituting the obtained
expressions for $x,y$ into the invariant derivative $\nabla_1$ and simplifying
modulo the obtained invariants (note that $\nabla_2$ is already $H$-invariant)
we get new invariant derivative
 $$
\nabla'_1=(u_xu_{yy}-u_yu_{xy})\mathcal{D}_x-(u_xu_{xy}-u_yu_{xx})\mathcal{D}_y.
 $$
Note that $\nabla'_1(I_0)=I_{2c}$ so the latter generator can be omitted.
The commutator of invariant derivations is
 $$
[\nabla'_1,\nabla_2] = -\frac{\nabla_2(I_{2c})}{I_{2c}}\nabla'_1+
\Bigl(\frac{\nabla'_1(I_{2c})}{I_{2c}}-2I_2'\Bigr)\nabla_2.
 $$
Denote this relation and the relation $\nabla_2(I_0)=0$ by $\mathcal{R}_1,\mathcal{R}_2$. Note that $\nabla'_1(I_0)$ is a
second order invariant, and application of $\nabla'_1,\nabla_2$ to it
and $I_2'$ gives 4 third order invariants. Further differentiation
gives 6 fourth order invariants, whence the syzyzy $\mathcal{R}_3$:
 \begin{multline*}
-I_{2c}\nabla_{2}(I_{3b})+\nabla'_{1}(I_{3c})+I'_{2}\nabla_{2}(I_{3d})\\
-\frac{1}{I_{2c}}\Bigl[12I_{2}'^{2}I_{2c}^{2}-10I'_{2}I_{2c}I_{3c}
+3I'_{2}I_{3d}^{2}+3I_{2c}^{2}I_{3a}-3I_{2c}I_{3b}I_{3d}+3I_{3c}^{2}\Bigr]=0,
 \end{multline*}
where $I_{3a}=\nabla'_{1}(I'_{2}), I_{3b}=\nabla_{2}(I'_{2}), I_{3c}=\nabla'_{1}(I_{2c})$ and $I_{3d}=\nabla_{2}(I_{2c})$. Therefore the algebra of differential $H$-invariants is
\begin{framed}[.5\textwidth]
 $$
\mathcal{A}_H=\la I_0,I'_2\,;\,\nabla'_1,\nabla_2\,|\,
\mathcal{R}_1,\mathcal{R}_2,\mathcal{R}_3\ra.
 $$
\end{framed}

\subsection{Affine symplectic group action: curves}\label{asga2}

Now we discuss the case of curves on the conformal symplectic plane,
with the same $H$ as in \S\ref{asga1}.
Consider the invariants from Section \ref{C:n=1} and note that $I_{k+2}=\nabla^kI_2$
are not translationally invariant. However using the elimination of parameters trick as above
we arrive to micro-local differential invariant and invariant derivation
 $$
I'_4=\sqrt[3]{y_2}(3y_2^{-2}y_4-5y_2^{-3}y_3^2) ,\quad
\nabla'=\frac1{\sqrt[3]{y_2}}\mathcal{D}_x.
 $$
In other words, these are invariants with respect to $\mathfrak{h}$ but not with respect to $H$.
Indeed, by the global Lie-Tresse theorem \cite{KL} we know that
the invariants should be rational.
To get generators we therefore pass to
 $$
I_4''=(I_4')^3\quad\text{ and }\quad\nabla''=I_4'\nabla'.
 $$
Consequently these freely generate the algebra of differential $H$-invariants:
\begin{framed}[.25\textwidth]
 $$
\mathcal{A}_H=\la I''_4\,;\,\nabla''\ra.
 $$
\end{framed}

\subsection{Affine conformal symplectic group action: functions}\label{s75}

Let us discuss differential invariants of functions on the affine
conformal symplectic plane, $H=\op{ACSp}(2n,\R)$.
For $n=1$ observe $H=\op{Aff}(2,\R)$.
We can combine the approaches of the previous two sections, for instance by
taking the affine symplectic differential invariants and restricting to those
of weight 0 with respect to the scaling by the center action.

Referring to the notations of \S\ref{asga1} we get that the weights
of $I_0,I'_2$ are $0,-4$, while that of $\nabla_1',\nabla_2$
are $-4,-2$ respectively.
Therefore the algebra of differential invariants $\mathcal{A}_H$
is generated by the invariant derivations
 $$
\nabla''_1=\frac1{I'_2}\nabla',\quad \nabla_2''=\frac{I_2'}{\nabla_2(I_2')}\nabla_2
 $$
and the differential invariants (derived invariants $\nabla_1''I_0$, $(\nabla_1'')^2I_0$, $\nabla_2''\nabla_1''I_0$ are omitted)
 $$
I_0,\quad  I_{3a}''=\frac{\nabla_1'(I_2')}{(I_2')^2},\quad I_{3b}''=\frac{(\nabla_2I_2')^2}{(I_2')^3}.
 $$
Denote $\mathcal{R}_{l}$ as the unknown differential syzygies, then the algebra of differential invariants is generated as
\begin{framed}[.45\textwidth]
$$\mathcal{A}_{H}=\la
I_0, I_{3a}'', I_{3b}''\, ;\, \nabla_{1}'',\nabla_{2}'' \ |\ \mathcal{R}_{l}\ra. $$
\end{framed}

\subsection{Affine conformal symplectic group action: curves}

Similarly for the case of curves on the conformal symplectic plane,
with the same $H$ as in \S\ref{s75}, we get in the notations of
\S\ref{asga2} that the weights of $I_4''$ is $-4$ and
that of $\nabla''$ is $-2$.
Therefore the algebra of differential invariants $\mathcal{A}_H$
is generated by
 $$
I_5=\frac{(\nabla''I_4'')^2}{(I_4'')^3}\quad\text{ and }\quad\nabla'''=\frac{I_4''}{\nabla''(I_4'')}\nabla''.
 $$
In fact, it is a free differential algebra
\begin{framed}[.25\textwidth]
$$\mathcal{A}_H=\la I_5\,;\,\nabla'''\ra.$$
\end{framed}

\section{Differential invariants in Contact Spaces}

Let $W$ be a contact space that is a contactification
of the symplectic vector space $V$. In coordinates,
$W=\R^{2n+1}(\bx,\by,z)$ is equipped with the contact form
$\alpha=dz-\by\,d\bx$ such that its differential $d\alpha=d\bx\wedge d\by$
descends to the symplectic form on $V=\R^{2n}(\bx,\by)$.

As the equivalence group we take either $G=\op{Sp}(2n,\R)$
lifted to an action on $W$ from the standard linear action on $V$,
or its central extension $\hat G=\op{CSp}(2n,\R)$
corresponding to the scaling
$(\bx,\by,z)\mapsto(\lambda\bx,\lambda\by,\lambda^2z)$.
(One can also consider the affine extensions, as was done in
Section \ref{S7} but we skip doing this.)

Note that the group $G$ does not have an open orbit on $W$ because
$I_0=2z-\bx\by$ is an invariant. This gives a way to carry over
the results on the algebra of differential invariants in $V$ to
that in $W$ (for both functions and submanifolds; note that the formulae
from the symplectic case enter through a change of variables,
which is due to the lift of Hamiltonian vector fields $X_H$ to
contact Hamiltonian fields).

Then we can single out the subalgebra $\mathcal{A}_{\hat G}\subset
\mathcal{A}_G$ as the space of functions of weight 0 with respect to the
scaling above (or its infinitesimal field).
In particular, as $I_0$ has weight 2, it is not a scaling invariant,
and in fact the action of $\hat G$ on $W$ is almost transitive.

Below we demonstrate this two-stage computation in the simplest case $n=1$.
Note that the action of $\hat G=\op{GL}(2,\R)\supset G=\op{SL}(2,\R)$
on $W=\R^3(x,y,z)$ has the formula
 $$
\Phi_A(x,y,z)=\bigl(ax+by,cx+dy,(ad-bc)(z-\tfrac12xy)
+\tfrac12(ax+by)(cx+dy)\bigr).
 $$
with $A=\begin{pmatrix}a&b\\ c&d\end{pmatrix}\in\hat G$.
This explicit parametrization is a base for an application
of the moving frame method, which involves normalization of the
group parameters via elimination. (This was already exploited in
\S\ref{asga1}-\ref{asga2}.) This algorithm
(we refer for details to \cite{Ol}; an elaborated version of it,
the method of equivariant moving frame, was further developed
in the works by Peter Olver and co-authors)
allows to carry the computations below,
however for $n>1$ it would meet the complexity issues.
Yet the method we propose works for arbitrary $n>1$
as a straightforward generalization.

\subsection{Differential Invariants: Curves}

We begin with the group $G=\op{Sp}(2,\R)=\op{SL}(2,\R)$.
Its action on $W$ has a base invariant $I_0=2z-xy$.

The curves will be represented as $y=y(x)$, $z=z(x)$ and the projection
to $\R^2(x,y)$ restores the symplectic action. We note that
invariants from \S\ref{C:n=1} are still $G$-invariants in the contact action,
and we will use them:
 $$
I_{2a}=\frac{y_2}{(xy_1-y)^3},\qquad \nabla=\frac1{xy_1-y}\mathcal{D}_x.
 $$
Differential invariants of order $\leq 2$ are generated by $I_0$, $I_1=\nabla(I_0)$,
$I_{2a}$ and $I_{2b}=\nabla(I_1)$. Of course, in Lie-Tresse generating set we
omit the derived invariants $I_1,I_{2b}$, namely
\begin{framed}[.3\textwidth]
 $$
\mathcal{A}_G=\la I_0,I_{2a}\,;\,\nabla\ra.
 $$
\end{framed}
However these derived invariants are useful in generating the algebra of $\hat G$-invariants.
Indeed, with respect to the action of the center $\xi=x\p_x+y\p_y+2z\p_z$,
the weights of $I_0,I_1,I_{2a},I_{2b}$ are $2,0,-4,-2$ and
the weight of $\nabla$ is $-2$.
Thus, in order to obtain $\hat G$-invariants we pass to weight 0 combinations
($I_1$ is already invariant)
 $$
I'_{2a}=I_0^2I_{2a},\ I'_{2b}=I_0I_{2b},\quad \nabla'=I_0\nabla.
 $$
Explicitly after simplifications $I_1\mapsto\frac12(I_1+1)$, $I_{2b}'\mapsto\frac12I_{2b}'$ we get:
 \begin{align*}
I_1&=\frac{z_1-y}{xy_1-y},\qquad \nabla' = \frac{2z-xy}{xy_1-y}\mathcal{D}_x,\\
I'_{2a}&=\frac{(2z-xy)^2}{(xy_1-y)^3}\,y_2,\\
I'_{2b}&=\frac{2z-xy}{(xy_1-y)^3}\Bigl(x(y-z_1)y_2-(xy_1-y)(y_1-z_2)\Bigr).
 \end{align*}

The count of invariants is $h_0=0$, $h_1=1$ and $h_k=2$ for $k\ge2$. We conclude:
 \begin{theorem}
The algebra of differential invariants of the $\hat G$-action on $J^\infty(W,1)$
is freely generated as follows:
 $$
\mathcal{A}_{\hat G}=\la I_1,I'_{2a}\,;\,\nabla'\ra.
 $$
 \end{theorem}

\subsection{Differential Invariants: Surfaces}

Now we consider the action of $G$ and $\hat G$ on surfaces given as $z=z(x,y)$.
Since projection to $\R^2(x,y)$ gives the symplectic plane the $G$-computations can
be derived from \S\ref{F:n=1} with substitution $u=2z-xy$. This gives us the following
differential invariants and invariant derivations with respect to $G$:
 $$
I_0=2z-xy,\qquad \nabla_1=x\mathcal{D}_x+y\mathcal{D}_y,\quad \nabla_2=(x-2z_y)\mathcal{D}_x+(2z_x-y)\mathcal{D}_y,
 $$
and with the notations $I_1=\frac12\nabla_1(I_0)$, $I_{2a}=\nabla_1(I_1)-I_1$, $I_{2b}=-\frac12(\nabla_2(I_1)+I_{2a}-I_1)$
the following first and second order invariants
 \begin{align*}
I_1 &= xz_x+yz_y-xy,\\
I_{2a} &= x^2z_{xx}+2xyz_{xy}+y^2z_{yy}-xy, \\
I_{2b} &= x(z_y-x)z_{xx}-yz_xz_{yy}+(y(z_y-x)-xz_x)z_{xy}+xz_x, \\
I_{2c} &= z_x^2z_{yy}-2z_x(z_y-x)z_{xy}+(z_y-x)^2z_{xx}+z_x(z_y-x).
 \end{align*}
Now to obtain $\hat G$-invariants note that $I_0,I_1,I_{2a},I_{2b},I_{2c}$ all have weight 2
with respect to $\xi$, while $\nabla_1,\nabla_2$ are already invariant.
Thus the invariants are
 $$
I'_1=I_0^{-1}I_1,\ I'_{2a}=I_0^{-1}I_{2a},\ I'_{2b}=I_0^{-1}I_{2b},\ I'_{2c}=I_0^{-1}I_{2c}.
 $$
We have:
 $$
I'_{2a}=\nabla_1(I'_1)+2(I'_1)^2-I'_1,\quad I'_{2b}=-\frac12\nabla_1(I'_1)-\frac12\nabla_2(I'_1)-(I_1')^2+I'_1,
 $$
so these can be omitted from the list of generators.

The count of $\hat G$-invariants is $h_0=0$, $h_1=1$ and $h_k=k+1$ for $k\ge2$.

Applying the derivations to the generating invariants and counting the relations,
we find that beside the commutation relation
 $$
[\nabla_1,\nabla_2]+\frac{\nabla_2(I'_1)}{I'_1}\nabla_1-\Bigl(\frac{\nabla_1(I'_1)}{I'_1}+2(I'_1-1)\Bigr)\nabla_2=0
 $$
there is one more relation generating the module of differential syzygies
 \begin{multline*}
\nabla_1^2(I'_1)+2\nabla_1\nabla_2(I'_1)+\nabla_2^2(I'_1)-4\nabla_1(I'_{2c})\\
-3(I'_1)^{-1}\bigl(\nabla_1(I'_1)^2+2\nabla_1(I'_1)\nabla_2(I'_1)-4\nabla_1(I'_1)I'_{2c}+\nabla_2(I'_1)^2\bigr)\\
-2(I'_1-1)\bigl(3\nabla_1(I'_1)+4\nabla_2(I'_1)-8I'_{2c}\bigr)-4I'_1(I'_1-1)(2I'_1-1)=0.
 \end{multline*}
Denote these syzygies by $\mathcal{R}_1$ and $\mathcal{R}_2$.

Let us summarize the results.

 \begin{theorem}
The algebra of differential invariants of the $\hat G$-action on $J^\infty(W,2)$
is generated as follows:
 $$
\mathcal{A}_{\hat G}=\la I'_1,I'_{2c}\,;\,\nabla_1,\nabla_2\ |\ \mathcal{R}_1,\mathcal{R}_2\ra.
 $$
 \end{theorem}

\subsection{Differential Invariants: Functions}

Skipping the intermediate computation with the group $G$ let us directly pass to
the description of invariants on $J^\infty(W)$ with respect to the group $\hat G$.
Fix the coordinates as follows: $W=\R^3(x,y,z)$ with the contact form
$\alpha=dz-y\,dx$ as before, $J^0=W\times \R(u)$ and for the jet-coordinates
we use the numbered multi-index notations $u_\sigma$.

The count of the number of differential invariants is as follows:
$h_0=1$, $h_1=2$ and $h_k=\binom{k+2}2$ for $k\ge2$.

The zero and first order invariants are
 $$
I_0=u,\quad I_{1a}=(xy-2z)u_3, \quad I_{1b}= xu_1+y(u_2+xu_3).
 $$
Next we obtain the invariant derivations
 \begin{align*}
\nabla_1 &= x\mathcal{D}_x+y\mathcal{D}_y+2z\mathcal{D}_z \\
\nabla_2 &= (xy-2z)\mathcal{D}_z \\
\nabla_3 &= (xu_3+u_2)(xy-2z)\mathcal{D}_x-u_1(xy-2z)\mathcal{D}_y-xu_1(xy-2z)\mathcal{D}_z
 \end{align*}
Note that $I_{1a}=\nabla_2(I_0)$, $I_{1b}=(\nabla_1+\nabla_2)(I_0)$ and $\nabla_3(I_0)=0$.
The latter is the first differential syzygy, denoted
 $$
\mathcal{R}_1=\nabla_3(I_0).
 $$

Passing to $J^2W$ we can produce 6 second order differential invariants $\nabla_i(I_{1a}),\nabla_i(I_{1b})$ but only 5 of those
are independent. We find the remaining 1 differential invariant via the the method of moving frames and get


 \begin{longtable}{|l|}
\hline
Second Order Differential Invariants  \\ \hline\\
$I_{2a}=y^2u_{2,2}+y(4zu_{2,3}+2xu_{1,2}+u_2)+4z^2u_{3,3}$ \\
$\hspace{0.55cm}+z(4xu_{1,3}+4u_3)+x(xu_{1,1}+u_1)$ \\
\empty \\
$I_{2b}=(xy-2z)(yu_{2,3}+2zu_{3,3}+xu_{1,3}+2u_3)$  \\
\empty \\
$I_{2c}=-(xy-2z)(x^2(u_1u_{1,3}-u_3u_{1,1})+x(u_1(yu_{2,3}+2zu_{3,3}+u_3+u_{1,2})$ \\
$\hspace{0.55cm}-yu_3u_{1,2}-2zu_3u_{1,3}-u_2u_{1,1})+u_1(yu_{2,2}+2zu_{2,3})$ \\
$\hspace{0.55cm}-u_2(yu_{1,2}+2zu_{1,3}))$ \\
\empty \\
$I_{2d}=(xy-2z)(-2u_3+(xy-2z)u_{3,3})$ \\
\empty \\
$I_{2e}=-(xy-2z)(x^2y(u_1u_{3,3}-u_3u_{1,3})+x(y(u_1u_{2,3}-u_3^2-u_2u_{1,3})$ \\
$\hspace{0.55cm}+u_1(-2zu_{3,3}-u_3)+2zu_3u_{1,3})-yu_2u_3-2z(u_1u_{2,3}-u_2u_{1,3}))$ \\
\empty \\
$I_{2f}=x^4y^2u_1^2u_{3,3}-2x^4y^2u_1u_3u_{1,3}+x^4y^2u_3^2u_{1,1}+2x^3y^2u_1^2u_{2,3}$ \\
$\hspace{0.55cm}-x^3y^2u_1u_3^2-2x^3y^2u_1u_3u_{1,2}-2x^3y^2u_1u_2u_{1,3}+2x^3y^2u_2u_3u_{1,1}$ \\
$\hspace{0.55cm}-4x^3yzu_1^2u_{3,3}+8x^3yzu_1u_3u_{1,3}-4x^3yzu_3^2u_{1,1}+x^2y^2u_1^2u_{2,2}$ \\
$\hspace{0.55cm}-x^2y^2u_1u_2u_3-2x^2y^2u_1u_2u_{1,2}+x^2y^2u_2^2u_{1,1}-8x^2yzu_1^2u_{2,3}$ \\
$\hspace{0.55cm}+4x^2yzu_1u_3^2+8x^2yzu_1u_3u_{1,2}+8x^2yzu_1u_2u_{1,3}-8x^2yzu_2u_3u_{1,1}$ \\
$\hspace{0.55cm}+4x^2z^2u_1^2u_{3,3}-8x^2z^2u_1u_3u_{1,3}+4x^2z^2u_3^2u_{1,1}-4xyzu_1^2u_{2,2}$ \\
$\hspace{0.55cm}+4xyzu_1u_2u_3+8xyzu_1u_2u_{1,2}-4xyzu_2^2u_{1,1}+8xz^2u_1^2u_{2,3}$ \\
$\hspace{0.55cm}-4xz^2u_1u_3^2-8xz^2u_1u_3u_{1,2}-8xz^2u_1u_2u_{1,3}+8xz^2u_2u_3u_{1,1}$ \\
$\hspace{0.55cm}+4z^2u_1^2u_{2,2}-4z^2u_1u_2u_3-8z^2u_1u_2u_{1,2}+4z^2u_2^2u_{1,1}$ \\ 
\empty \\
\hline
 \end{longtable}

Note that $I_{2a},I_{2b},I_{2c},I_{2d},I_{2e}$ can be expressed through $I_0,I_{1a},I_{1b}$ and invariant derivations.
Thus they need not enter the set of generators.

All the differential syzygies coming from the commutators are
\begin{align*}
    \mathcal{R}_{2}&=[\nabla_{1},\nabla_{2}] \\
    \mathcal{R}_{3}&=(I_{1a}+I_{1b})[\nabla_{1},\nabla_{3}]+I_{2c}(\nabla_{1}+\nabla_{2})-(I_{2a}+I_{2b})\nabla_{3} \\
    \mathcal{R}_{4}&=(I_{1a}+I_{1b})[\nabla_{2},\nabla_{3}]-(I_{1b}(I_{1a}+I_{1b})-I_{2e})\nabla_{1}+(I_{1a}(I_{1a}+I_{1b})+I_{2e})\nabla_{2}\\
    &-(I_{2b}+I_{2d}-2(I_{1a}+I_{1b}))\nabla_{3}
\end{align*}
The remaining differential syzygies are found by the symbolic method: find a relation between the symbols of differentiated invariants,
get a linear combination of lower order and express it through the invariants established earlier.
 \begin{align*}
\mathcal{R}_{5}&=(I_{1a}+I_{1b})(\nabla_{3}(I_{2b})-\nabla_{1}(I_{2e}))-(I_{2c}-I_{2e})I_{2b}+I_{2a}I_{2e}-I_{2c}I_{2d} \\
\mathcal{R}_{6}&=(I_{1a}+I_{1b})(\nabla_{3}(I_{2c})-\nabla_{1}(I_{2f}))-3I_{2c}^{2}-(I_{1a}^{2}+I_{1a}I_{1b}+3I_{2e})I_{2c}+3I_{2f}(I_{2a}+I_{2b}) \\
\mathcal{R}_{7}&=(I_{1a}+I_{1b})(-\nabla_{3}(I_{2e})+\nabla_{2}(I_{2f}))-I_{1b}^{4}-4I_{1a}I_{1b}^{3}-(5I_{1a}^{2}+2I_{2c})I_{1b}^{2}\\
    &-(2I_{1a}^{3}+(2I_{2c}-3I_{2e})I_{1a}+4I_{2f})I_{1b}+3I_{2e}I_{1a}^{2}+4I_{2f}I_{1a}+3I_{2e}^{2}\\
    &+3I_{2c}I_{2e}-3I_{2f}(I_{2b}+I_{2d})
 \end{align*}

Let us summarize the results.

 \begin{theorem}
The algebra of differential invariants of the $\hat G$-action on $J^\infty(W)$
is generated as follows:
 $$
\mathcal{A}_{\hat G}=\la I_0,I_{2f} \,;\, \nabla_1,\nabla_2,\nabla_3 \ |\ \mathcal{R}_i=0,\ i=1\hdots7 \ra.
 $$
 \end{theorem}

\section{Conclusion}

In this paper we computed the algebra of differential invariants
for various geometric objects on symplectic spaces with several choices
of the equivalence group and touched upon a relation between the
invariants of the pair (group,subgroup) action.

For most of the text we worked with the linear symplectic group, but
we demonstrated how to extend the results for conformal symplectic
and affine symplectic groups, treated in other publications.
Some of the objects were also investigated by different authors, namely
jets of curves \cite{KOT,V} and hypersurfaces \cite{D}, yet the technique
and the description of the algebras are quite distinct.
Surfaces in four-dimensional symplectic space were also studied
in \cite{MK,MN,MH}, but they considered Lagrangian surfaces while
our focus was on symplectic (generic) submanifolds.

Other geometric objects appeared in \cite{X}, which intersects with our
work by studying functions on the symplectic spaces.
Again the approaches differ significantly:
in \cite{X} the infinite number of generators were computed
(with a nontrivial change of variables)
while our method uses the Lie-Tresse finite type presentation
of the algebra (in the original jet-coordinates).
This latter allows, in particular, to solve the equivalence problem
via a finite-dimensional signature variety.

The work \cite{X} also described invariants in the adjoint bundle,
and one can consider other geometric spaces on which the symplectic group
acts. For instance, \cite{DFKN} was devoted to four-fold surfaces
in 6-dimensional Lagrangian Grassmanian, satisfying the
integrability condition.
It would be worth characterizing those via symplectic invariants.

Finally note that one can approach the equivalence problem
of geometric objects via discretizations, with more algebraic
methods, see \cite{AK}.

\appendix
\section{Differential Invariants of curves in 4-dimensions}

Here are explicit expressions of the differential invariants of curves $x=x(t),y=y(t),z=z(t)$, as derived in Section \ref{C:n=2}.
These as well as other long formulae resulting from our calculations can be also found in ancillary files of this arXiv submission.

Below $\gamma=1/(ty_{1}+xz_{1}-x_{1}z-y)$ is the factor of $\mathcal{D}_t$ in $\nabla$. The jet notations are $x_{t}=x_{1},x_{tt}=x_{2},x_{ttt}=x_{3}$ etc, likewise for $y$ and $z$.


\begin{longtable}{|l|}
\hline
Differential invariants that together with $\nabla$ generate $\mathcal{A}$ \\ \hline\\
$I_{2}=\gamma^{3}(x_{1}z_{2}-z_{1}x_{2}+y_{2})$ \\
\empty \\
$I_{3b}=-\gamma^{6}(tx_{1}y_{2}z_{3}-tx_{1}y_{3}z_{2}-tx_{2}y_{1}z_{3}+tx_{2}y_{3}z_{1}+tx_{3}y_{1}z_{2}-tx_{3}y_{2}z_{1}$ \\
$\hspace{0.55cm}-xy_{2}z_{3}+xy_{3}z_{2}+x_{2}yz_{3}-x_{2}y_{3}z-x_{3}yz_{2}+x_{3}y_{2}z)$ \\
\empty \\
$I_{4c}=-\gamma^{10}(t^{3}x_{1}y_{1}^{2}y_{3}z_{4}-t^{3}x_{1}y_{1}^{2}y_{4}z_{3}-3t^{3}x_{1}y_{1}y_{2}^{2}z_{4}+4t^{3}x_{1}y_{1}y_{2}y_{3}z_{3}$ \\
$\hspace{0.55cm}+3t^{3}x_{1}y_{1}y_{2}y_{4}z_{2}-4t^{3}x_{1}y_{1}y_{3}^{2}z_{2}+3t^{3}x_{1}y_{2}^{3}z_{3}-3t^{3}x_{1}y_{2}^{2}y_{3}z_{2}+3t^{3}x_{2}y_{1}^{2}y_{2}z_{4}$ \\
$\hspace{0.55cm}-4t^{3}x_{2}y_{1}^{2}y_{3}z_{3}-3t^{3}x_{2}y_{1}y_{2}^{2}z_{3}-3t^{3}x_{2}y_{1}y_{2}y_{4}z_{1}+4t^{3}x_{2}y_{1}y_{3}^{2}z_{1}$ \\
$\hspace{0.55cm}+3t^{3}x_{2}y_{2}^{2}y_{3}z_{1}-t^{3}x_{3}y_{1}^{3}z_{4}+4t^{3}x_{3}y_{1}^{2}y_{3}z_{2}+t^{3}x_{3}y_{1}^{2}y_{4}z_{1}+3t^{3}x_{3}y_{1}y_{2}^{2}z_{2}$ \\
$\hspace{0.55cm}-4t^{3}x_{3}y_{1}y_{2}y_{3}z_{1}-3t^{3}x_{3}y_{2}^{3}z_{1}+t^{3}x_{4}y_{1}^{3}z_{3}-3t^{3}x_{4}y_{1}^{2}y_{2}z_{2}-t^{3}x_{4}y_{1}^{2}y_{3}z_{1}$ \\
$\hspace{0.55cm}+3t^{3}x_{4}y_{1}y_{2}^{2}z_{1}-3t^{2}xx_{1}y_{1}y_{2}z_{2}z_{4}+4t^{2}xx_{1}y_{1}y_{2}z_{3}^{2}+2t^{2}xx_{1}y_{1}y_{3}z_{1}z_{4}$ \\
$\hspace{0.55cm}-4t^{2}xx_{1}y_{1}y_{3}z_{2}z_{3}-2t^{2}xx_{1}y_{1}y_{4}z_{1}z_{3}+3t^{2}xx_{1}y_{1}y_{4}z_{2}^{2}
-3t^{2}xx_{1}y_{2}^{2}z_{1}z_{4}$ \\
$\hspace{0.55cm}+6t^{2}xx_{1}y_{2}^{2}z_{2}z_{3}+4t^{2}xx_{1}y_{2}y_{3}z_{1}z_{3}-6t^{2}xx_{1}y_{2}y_{3}z_{2}^{2}+3t^{2}xx_{1}y_{2}y_{4}z_{1}z_{2}$ \\
$\hspace{0.55cm}-4t^{2}xx_{1}y_{3}^{2}z_{1}z_{2}+3t^{2}xx_{2}y_{1}^{2}z_{2}z_{4}-4t^{2}xx_{2}y_{1}^{2}z_{3}^{2}+3t^{2}xx_{2}y_{1}y_{2}z_{1}z_{4}$ \\
$\hspace{0.55cm}-6t^{2}xx_{2}y_{1}y_{2}z_{2}z_{3}-3t^{2}xx_{2}y_{1}y_{4}z_{1}z_{2}+6t^{2}xx_{2}y_{2}y_{3}z_{1}z_{2}-3t^{2}xx_{2}y_{2}y_{4}z_{1}^{2}$ \\
$\hspace{0.55cm}+4t^{2}xx_{2}y_{3}^{2}z_{1}^{2}-2t^{2}xx_{3}y_{1}^{2}z_{1}z_{4}+4t^{2}xx_{3}y_{1}^{2}z_{2}z_{3}-4t^{2}xx_{3}y_{1}y_{2}z_{1}z_{3}$ \\
$\hspace{0.55cm}+6t^{2}xx_{3}y_{1}y_{2}z_{2}^{2}+4t^{2}xx_{3}y_{1}y_{3}z_{1}z_{2}+2t^{2}xx_{3}y_{1}y_{4}z_{1}^{2}-6t^{2}xx_{3}y_{2}^{2}z_{1}z_{2}$ \\
$\hspace{0.55cm}-4t^{2}xx_{3}y_{2}y_{3}z_{1}^{2}+2t^{2}xx_{4}y_{1}^{2}z_{1}z_{3}-3t^{2}xx_{4}y_{1}^{2}z_{2}^{2}-2t^{2}xx_{4}y_{1}y_{3}z_{1}^{2}$ \\
$\hspace{0.55cm}+3t^{2}xx_{4}y_{2}^{2}z_{1}^{2}-2t^{2}x_{1}^{2}y_{1}y_{3}zz_{4}+2t^{2}x_{1}^{2}y_{1}y_{4}zz_{3}+3t^{2}x_{1}^{2}y_{2}^{2}zz_{4}$ \\
$\hspace{0.55cm}-4t^{2}x_{1}^{2}y_{2}y_{3}zz_{3}-3t^{2}x_{1}^{2}y_{2}y_{4}zz_{2}+4t^{2}x_{1}^{2}y_{3}^{2}zz_{2}+4t^{2}x_{1}x_{2}y_{1}y_{3}zz_{3}$ \\
$\hspace{0.55cm}-3t^{2}x_{1}x_{2}y_{1}y_{4}zz_{2}-6t^{2}x_{1}x_{2}y_{2}^{2}zz_{3}+6t^{2}x_{1}x_{2}y_{2}y_{3}zz_{2}+3t^{2}x_{1}x_{2}y_{2}y_{4}zz_{1}$ \\
$\hspace{0.55cm}-4t^{2}x_{1}x_{2}y_{3}^{2}zz_{1}+2t^{2}x_{1}x_{3}y_{1}^{2}zz_{4}-4t^{2}x_{1}x_{3}y_{1}y_{2}zz_{3}-2t^{2}x_{1}x_{3}y_{1}y_{4}zz_{1}$ \\
$\hspace{0.55cm}+4t^{2}x_{1}x_{3}y_{2}y_{3}zz_{1}-2t^{2}x_{1}x_{4}y_{1}^{2}zz_{3}+3t^{2}x_{1}x_{4}y_{1}y_{2}zz_{2}+2t^{2}x_{1}x_{4}y_{1}y_{3}zz_{1}$ \\
$\hspace{0.55cm}-3t^{2}x_{1}x_{4}y_{2}^{2}zz_{1}-3t^{2}x_{2}^{2}y_{1}^{2}zz_{4}+6t^{2}x_{2}^{2}y_{1}y_{2}zz_{3}+3t^{2}x_{2}^{2}y_{1}y_{4}zz_{1}$ \\
$\hspace{0.55cm}-6t^{2}x_{2}^{2}y_{2}y_{3}zz_{1}+4t^{2}x_{2}x_{3}y_{1}^{2}zz_{3}-6t^{2}x_{2}x_{3}y_{1}y_{2}zz_{2}-4t^{2}x_{2}x_{3}y_{1}y_{3}zz_{1}$ \\
$\hspace{0.55cm}+6t^{2}x_{2}x_{3}y_{2}^{2}zz_{1}+3t^{2}x_{2}x_{4}y_{1}^{2}zz_{2}-3t^{2}x_{2}x_{4}y_{1}y_{2}zz_{1}-4t^{2}x_{3}^{2}y_{1}^{2}zz_{2}$ \\
$\hspace{0.55cm}+4t^{2}x_{3}^{2}y_{1}y_{2}zz_{1}-3tx^{2}x_{1}y_{2}z_{1}z_{2}z_{4}+4tx^{2}x_{1}y_{2}z_{1}z_{3}^{2}+3tx^{2}x_{1}y_{2}z_{2}^{2}z_{3}$ \\
$\hspace{0.55cm}+tx^{2}x_{1}y_{3}z_{1}^{2}z_{4}-4tx^{2}x_{1}y_{3}z_{1}z_{2}z_{3}-3tx^{2}x_{1}y_{3}z_{2}^{3}-tx^{2}x_{1}y_{4}z_{1}^{2}z_{3}$ \\
$\hspace{0.55cm}+3tx^{2}x_{1}y_{4}z_{1}z_{2}^{2}+3tx^{2}x_{2}y_{1}z_{1}z_{2}z_{4}-4tx^{2}x_{2}y_{1}z_{1}z_{3}^{2}-3tx^{2}x_{2}y_{1}z_{2}^{2}z_{3}$ \\
$\hspace{0.55cm}+4tx^{2}x_{2}y_{3}z_{1}^{2}z_{3}+3tx^{2}x_{2}y_{3}z_{1}z_{2}^{2}-3tx^{2}x_{2}y_{4}z_{1}^{2}z_{2}-tx^{2}x_{3}y_{1}z_{1}^{2}z_{4}$ \\
$\hspace{0.55cm}+4tx^{2}x_{3}y_{1}z_{1}z_{2}z_{3}+3tx^{2}x_{3}y_{1}z_{2}^{3}-4tx^{2}x_{3}y_{2}z_{1}^{2}z_{3}-3tx^{2}x_{3}y_{2}z_{1}z_{2}^{2}+tx^{2}x_{3}y_{4}z_{1}^{3}$ \\
$\hspace{0.55cm}+tx^{2}x_{4}y_{1}z_{1}^{2}z_{3}-3tx^{2}x_{4}y_{1}z_{1}z_{2}^{2}+3tx^{2}x_{4}y_{2}z_{1}^{2}z_{2}-tx^{2}x_{4}y_{3}z_{1}^{3}+3txx_{1}^{2}y_{2}zz_{2}z_{4}$ \\
$\hspace{0.55cm}-4txx_{1}^{2}y_{2}zz_{3}^{2}-2txx_{1}^{2}y_{3}zz_{1}z_{4}+4txx_{1}^{2}y_{3}zz_{2}z_{3}+2txx_{1}^{2}y_{4}zz_{1}z_{3}-3txx_{1}^{2}y_{4}zz_{2}^{2}$ \\
$\hspace{0.55cm}-3txx_{1}x_{2}y_{1}zz_{2}z_{4}+4txx_{1}x_{2}y_{1}zz_{3}^{2}+3txx_{1}x_{2}y_{2}zz_{1}z_{4}-6txx_{1}x_{2}y_{2}zz_{2}z_{3}$ \\
$\hspace{0.55cm}-4txx_{1}x_{2}y_{3}zz_{1}z_{3}+6txx_{1}x_{2}y_{3}zz_{2}^{2}+2txx_{1}x_{3}y_{1}zz_{1}z_{4}-4txx_{1}x_{3}y_{1}zz_{2}z_{3}$ \\
$\hspace{0.55cm+4txx_{1}x_{3}y_{3}zz_{1}z_{2}-2txx_{1}x_{3}y_{4}zz_{1}^{2}}-2txx_{1}x_{4}y_{1}zz_{1}z_{3}+3txx_{1}x_{4}y_{1}zz_{2}^{2}$ \\
$\hspace{0.55cm}-3txx_{1}x_{4}y_{2}zz_{1}z_{2}+2txx_{1}x_{4}y_{3}zz_{1}^{2}-3txx_{2}^{2}y_{1}zz_{1}z_{4}+6txx_{2}^{2}y_{1}zz_{2}z_{3}$ \\
$\hspace{0.55cm}-6txx_{2}^{2}y_{3}zz_{1}z_{2}+3txx_{2}^{2}y_{4}zz_{1}^{2}+4txx_{2}x_{3}y_{1}zz_{1}z_{3}-6txx_{2}x_{3}y_{1}zz_{2}^{2}$ \\
$\hspace{0.55cm}+6txx_{2}x_{3}y_{2}zz_{1}z_{2}-4txx_{2}x_{3}y_{3}zz_{1}^{2}+3txx_{2}x_{4}y_{1}zz_{1}z_{2}-3txx_{2}x_{4}y_{2}zz_{1}^{2}$ \\
$\hspace{0.55cm}-4txx_{3}^{2}y_{1}zz_{1}z_{2}+4txx_{3}^{2}y_{2}zz_{1}^{2}+tx_{1}^{3}y_{3}z^{2}z_{4}-tx_{1}^{3}y_{4}z^{2}z_{3}-3tx_{1}^{2}x_{2}y_{2}z^{2}z_{4}$ \\
$\hspace{0.55cm}+3tx_{1}^{2}x_{2}y_{4}z^{2}z_{2}-tx_{1}^{2}x_{3}y_{1}z^{2}z_{4}+4tx_{1}^{2}x_{3}y_{2}z^{2}z_{3}-4tx_{1}^{2}x_{3}y_{3}z^{2}z_{2}+tx_{1}^{2}x_{3}y_{4}z^{2}z_{1}$ \\
$\hspace{0.55cm}+tx_{1}^{2}x_{4}y_{1}z^{2}z_{3}-tx_{1}^{2}x_{4}y_{3}z^{2}z_{1}+3tx_{1}x_{2}^{2}y_{1}z^{2}z_{4}+3tx_{1}x_{2}^{2}y_{2}z^{2}z_{3}-3tx_{1}x_{2}^{2}y_{3}z^{2}z_{2}$ \\
$\hspace{0.55cm}-3tx_{1}x_{2}^{2}y_{4}z^{2}z_{1}-4tx_{1}x_{2}x_{3}y_{1}z^{2}z_{3}+4tx_{1}x_{2}x_{3}y_{3}z^{2}z_{1}-3tx_{1}x_{2}x_{4}y_{1}z^{2}z_{2}$ \\
$\hspace{0.55cm}+3tx_{1}x_{2}x_{4}y_{2}z^{2}z_{1}+4tx_{1}x_{3}^{2}y_{1}z^{2}z_{2}-4tx_{1}x_{3}^{2}y_{2}z^{2}z_{1}-3tx_{2}^{3}y_{1}z^{2}z_{3}+3tx_{2}^{3}y_{3}z^{2}z_{1}$ \\
$\hspace{0.55cm}+3tx_{2}^{2}x_{3}y_{1}z^{2}z_{2}-3tx_{2}^{2}x_{3}y_{2}z^{2}z_{1}-t^{2}xy_{1}^{2}y_{3}z_{4}+t^{2}xy_{1}^{2}y_{4}z_{3}+3t^{2}xy_{1}y_{2}^{2}z_{4}$ \\
$\hspace{0.55cm}-4t^{2}xy_{1}y_{2}y_{3}z_{3}-3t^{2}xy_{1}y_{2}y_{4}z_{2}+4t^{2}xy_{1}y_{3}^{2}z_{2}-3t^{2}xy_{2}^{3}z_{3}+3t^{2}xy_{2}^{2}y_{3}z_{2}$ \\
$\hspace{0.55cm}-2t^{2}x_{1}yy_{1}y_{3}z_{4}+2t^{2}x_{1}yy_{1}y_{4}z_{3}+3t^{2}x_{1}yy_{2}^{2}z_{4}-4t^{2}x_{1}yy_{2}y_{3}z_{3}-3t^{2}x_{1}yy_{2}y_{4}z_{2}$ \\
$\hspace{0.55cm}+4t^{2}x_{1}yy_{3}^{2}z_{2}-6t^{2}x_{2}yy_{1}y_{2}z_{4}+8t^{2}x_{2}yy_{1}y_{3}z_{3}+3t^{2}x_{2}yy_{2}^{2}z_{3}+3t^{2}x_{2}yy_{2}y_{4}z_{1}$ \\
$\hspace{0.55cm}-4t^{2}x_{2}yy_{3}^{2}z_{1}+3t^{2}x_{2}y_{1}y_{2}y_{4}z-4t^{2}x_{2}y_{1}y_{3}^{2}z-3t^{2}x_{2}y_{2}^{2}y_{3}z+3t^{2}x_{3}yy_{1}^{2}z_{4}$ \\
$\hspace{0.55cm}-8t^{2}x_{3}yy_{1}y_{3}z_{2}-2t^{2}x_{3}yy_{1}y_{4}z_{1}-3t^{2}x_{3}yy_{2}^{2}z_{2}+4t^{2}x_{3}yy_{2}y_{3}z_{1}$ \\
$\hspace{0.55cm}-t^{2}x_{3}y_{1}^{2}y_{4}z+4t^{2}x_{3}y_{1}y_{2}y_{3}z+3t^{2}x_{3}y_{2}^{3}z-3t^{2}x_{4}yy_{1}^{2}z_{3}+6t^{2}x_{4}yy_{1}y_{2}z_{2}$ \\
$\hspace{0.55cm}+2t^{2}x_{4}yy_{1}y_{3}z_{1}-3t^{2}x_{4}yy_{2}^{2}z_{1}+t^{2}x_{4}y_{1}^{2}y_{3}z-3t^{2}x_{4}y_{1}y_{2}^{2}z+3tx^{2}y_{1}y_{2}z_{2}z_{4}$ \\
$\hspace{0.55cm}-4tx^{2}y_{1}y_{2}z_{3}^{2}-2tx^{2}y_{1}y_{3}z_{1}z_{4}+4tx^{2}y_{1}y_{3}z_{2}z_{3}+2tx^{2}y_{1}y_{4}z_{1}z_{3}-3tx^{2}y_{1}y_{4}z_{2}^{2}$ \\
$\hspace{0.55cm}+3tx^{2}y_{2}^{2}z_{1}z_{4}-6tx^{2}y_{2}^{2}z_{2}z_{3}-4tx^{2}y_{2}y_{3}z_{1}z_{3}+6tx^{2}y_{2}y_{3}z_{2}^{2}-3tx^{2}y_{2}y_{4}z_{1}z_{2}$ \\
$\hspace{0.55cm}+4tx^{2}y_{3}^{2}z_{1}z_{2}+3txx_{1}yy_{2}z_{2}z_{4}-4txx_{1}yy_{2}z_{3}^{2}-2txx_{1}yy_{3}z_{1}z_{4}+4txx_{1}yy_{3}z_{2}z_{3}$ \\
$\hspace{0.55cm}+2txx_{1}yy_{4}z_{1}z_{3}-3txx_{1}yy_{4}z_{2}^{2}+2txx_{1}y_{1}y_{3}zz_{4}-2txx_{1}y_{1}y_{4}zz_{3}-3txx_{1}y_{2}^{2}zz_{4}$ \\
$\hspace{0.55cm}+4txx_{1}y_{2}y_{3}zz_{3}+3txx_{1}y_{2}y_{4}zz_{2}-4txx_{1}y_{3}^{2}zz_{2}-6txx_{2}yy_{1}z_{2}z_{4}+8txx_{2}yy_{1}z_{3}^{2}$ \\
$\hspace{0.55cm}-3txx_{2}yy_{2}z_{1}z_{4}+6txx_{2}yy_{2}z_{2}z_{3}+3txx_{2}yy_{4}z_{1}z_{2}-3txx_{2}y_{1}y_{2}zz_{4}-4txx_{2}y_{1}y_{3}zz_{3}$ \\
$\hspace{0.55cm}+6txx_{2}y_{1}y_{4}z_{{}}z_{2}+6txx_{2}y_{2}^{2}zz_{3}-12txx_{2}y_{2}y_{3}zz_{2}+3txx_{2}y_{2}y_{4}zz_{1}-4txx_{2}y_{3}^{2}zz_{1}$ \\
$\hspace{0.55cm}+4txx_{3}yy_{1}z_{1}z_{4}-8txx_{3}yy_{1}z_{2}z_{3}+4txx_{3}yy_{2}z_{1}z_{3}-6txx_{3}yy_{2}z_{2}^{2}-4txx_{3}yy_{3}z_{1}z_{2}$ \\
$\hspace{0.55cm}-2txx_{3}yy_{4}z_{1}^{2}+8txx_{3}y_{1}y_{2}zz_{3}-4txx_{3}y_{1}y_{3}zz_{2}-2txx_{3}y_{1}y_{4}zz_{1}+6txx_{3}y_{2}^{2}zz_{2}$ \\
$\hspace{0.55cm}+4txx_{3}y_{2}y_{3}zz_{1}-4txx_{4}yy_{1}z_{1}z_{3}+6txx_{4}yy_{1}z_{2}^{2}+2txx_{4}yy_{3}z_{1}^{2}-3txx_{4}y_{1}y_{2}zz_{2}$ \\
$\hspace{0.55cm}+2txx_{4}y_{1}y_{3}zz_{1}-3txx_{4}y_{2}^{2}zz_{1}+2tx_{1}^{2}yy_{3}zz_{4}-2tx_{1}^{2}yy_{4}zz_{3}-4tx_{1}x_{2}yy_{3}zz_{3}$ \\
$\hspace{0.55cm}+3tx_{1}x_{2}yy_{4}zz_{2}-3tx_{1}x_{2}y_{2}y_{4}z^{2}+4tx_{1}x_{2}y_{3}^{2}z^{2}-4tx_{1}x_{3}yy_{1}zz_{4}+4tx_{1}x_{3}yy_{2}zz_{3}$ \\
$\hspace{0.55cm}+2tx_{1}x_{3}yy_{4}zz_{1}+2tx_{1}x_{3}y_{1}y_{4}z^{2}-4tx_{1}x_{3}y_{2}y_{3}z^{2}+4tx_{1}x_{4}yy_{1}zz_{3}-3tx_{1}x_{4}yy_{2}zz_{2}$ \\
$\hspace{0.55cm}-2tx_{1}x_{4}yy_{3}zz_{1}-2tx_{1}x_{4}y_{1}y_{3}z^{2}+3tx_{1}x_{4}y_{2}^{2}z^{2}+6tx_{2}^{2}yy_{1}zz_{4}-6tx_{2}^{2}yy_{2}zz_{3}$ \\
$\hspace{0.55cm}-3tx_{2}^{2}yy_{4}zz_{1}-3tx_{2}^{2}y_{1}y_{4}z^{2}+6tx_{2}^{2}y_{2}y_{3}z^{2}-8tx_{2}x_{3}yy_{1}zz_{3}+6tx_{2}x_{3}yy_{2}zz_{2}$ \\
$\hspace{0.55cm}+4tx_{2}x_{3}yy_{3}zz_{1}+4tx_{2}x_{3}y_{1}y_{3}z^{2}-6tx_{2}x_{3}y_{2}^{2}z^{2}-6tx_{2}x_{4}yy_{1}zz_{2}+3tx_{2}x_{4}yy_{2}zz_{1}$ \\
$\hspace{0.55cm}+3tx_{2}x_{4}y_{1}y_{2}z^{2}+8tx_{3}^{2}yy_{1}zz_{2}-4tx_{3}^{2}yy_{2}zz_{1}-4tx_{3}^{2}y_{1}y_{2}z^{2}+3x^{3}y_{2}z_{1}z_{2}z_{4}$ \\
$\hspace{0.55cm}-4x^{3}y_{2}z_{1}z_{3}^{2}-3x^{3}y_{2}z_{2}^{2}z_{3}-x^{3}y_{3}z_{1}^{2}z_{4}+4x^{3}y_{3}z_{1}z_{2}z_{3}+3x^{3}y_{3}z_{2}^{3}+x^{3}y_{4}z_{1}^{2}z_{3}$ \\
$\hspace{0.55cm}-3x^{3}y_{4}z_{1}z_{2}^{2}-3x^{2}x_{1}y_{2}zz_{2}z_{4}+4x^{2}x_{1}y_{2}zz_{3}^{2}+2x^{2}x_{1}y_{3}zz_{1}z_{4}-4x^{2}x_{1}y_{3}zz_{2}z_{3}$ \\
$\hspace{0.55cm}-2x^{2}x_{1}y_{4}zz_{1}z_{3}+3x^{2}x_{1}y_{4}zz_{2}^{2}-3x^{2}x_{2}yz_{1}z_{2}z_{4}+4x^{2}x_{2}yz_{1}z_{3}^{2}+3x^{2}x_{2}yz_{2}^{2}z_{3}$ \\
$\hspace{0.55cm}-3x^{2}x_{2}y_{2}zz_{1}z_{4}+6x^{2}x_{2}y_{2}zz_{2}z_{3}-4x^{2}x_{2}y_{3}zz_{1}z_{3}-9x^{2}x_{2}y_{3}zz_{2}^{2}+6x^{2}x_{2}y_{4}zz_{1}z_{2}$ \\
$\hspace{0.55cm}+x^{2}x_{3}yz_{1}^{2}z_{4}-4x^{2}x_{3}yz_{1}z_{2}z_{3}-3x^{2}x_{3}yz_{2}^{3}+8x^{2}x_{3}y_{2}zz_{1}z_{3}+3x^{2}x_{3}y_{2}zz_{2}^{2}$ \\
$\hspace{0.55cm}-4x^{2}x_{3}y_{3}zz_{1}z_{2}-x^{2}x_{3}y_{4}zz_{1}^{2}-x^{2}x_{4}yz_{1}^{2}z_{3}+3x^{2}x_{4}yz_{1}z_{2}^{2}-3x^{2}x_{4}y_{2}zz_{1}z_{2}$ \\
$\hspace{0.55cm}+x^{2}x_{4}y_{3}zz_{1}^{2}-xx_{1}^{2}y_{3}z^{2}z_{4}+xx_{1}^{2}y_{4}z^{2}z_{3}+3xx_{1}x_{2}yzz_{2}z_{4}-4xx_{1}x_{2}yzz_{3}^{2}$ \\
$\hspace{0.55cm}+3xx_{1}x_{2}y_{2}z^{2}z_{4}+4xx_{1}x_{2}y_{3}z^{2}z_{3}-6xx_{1}x_{2}y_{4}z^{2}z_{2}-2xx_{1}x_{3}yzz_{1}z_{4}$ \\
$\hspace{0.55cm}+4xx_{1}x_{3}yzz_{2}z_{3}-8xx_{1}x_{3}y_{2}z^{2}z_{3}+4xx_{1}x_{3}y_{3}z^{2}z_{2}+2xx_{1}x_{3}y_{4}z^{2}z_{1}$ \\
$\hspace{0.55cm}+2xx_{1}x_{4}yzz_{1}z_{3}-3xx_{1}x_{4}yzz_{2}^{2}+3xx_{1}x_{4}y_{2}z^{2}z_{2}-2xx_{1}x_{4}y_{3}z^{2}z_{1}$ \\
$\hspace{0.55cm}+3xx_{2}^{2}yzz_{1}z_{4}-6xx_{2}^{2}yzz_{2}z_{3}-3xx_{2}^{2}y_{2}z^{2}z_{3}+9xx_{2}^{2}y_{3}z^{2}z_{2}-3xx_{2}^{2}y_{4}z^{2}z_{1}$ \\
$\hspace{0.55cm}-4xx_{2}x_{3}yzz_{1}z_{3}+6xx_{2}x_{3}yzz_{2}^{2}-6xx_{2}x_{3}y_{2}z^{2}z_{2}+4xx_{2}x_{3}y_{3}z^{2}z_{1}$ \\
$\hspace{0.55cm}-3xx_{2}x_{4}yzz_{1}z_{2}+3xx_{2}x_{4}y_{2}z^{2}z_{1}+4xx_{3}^{2}yzz_{1}z_{2}-4xx_{3}^{2}y_{2}z^{2}z_{1}+x_{1}^{2}x_{3}yz^{2}z_{4}$ \\
$\hspace{0.55cm}-x_{1}^{2}x_{3}y_{4}z^{3}-x_{1}^{2}x_{4}yz^{2}z_{3}+x_{1}^{2}x_{4}y_{3}z^{3}-3x_{1}x_{2}^{2}yz^{2}z_{4}+3x_{1}x_{2}^{2}y_{4}z^{3}$ \\
$\hspace{0.55cm}+4x_{1}x_{2}x_{3}yz^{2}z_{3}-4x_{1}x_{2}x_{3}y_{3}z^{3}+3x_{1}x_{2}x_{4}yz^{2}z_{2}-3x_{1}x_{2}x_{4}y_{2}z^{3}-4x_{1}x_{3}^{2}yz^{2}z_{2}$ \\
$\hspace{0.55cm}+4x_{1}x_{3}^{2}y_{2}z^{3}+3x_{2}^{3}yz^{2}z_{3}-3x_{2}^{3}y_{3}z^{3}-3x_{2}^{2}x_{3}yz^{2}z_{2}+3x_{2}^{2}x_{3}y_{2}z^{3}$ \\
$\hspace{0.55cm}+2txyy_{1}y_{3}z_{4}-2txyy_{1}y_{4}z_{3}-3txyy_{2}^{2}z_{4}+4txyy_{2}y_{3}z_{3}+3txyy_{2}y_{4}z_{2}$ \\
$\hspace{0.55cm}-4txyy_{3}^{2}z_{2}+tx_{1}y^{2}y_{3}z_{4}-tx_{1}y^{2}y_{4}z_{3}+3tx_{2}y^{2}y_{2}z_{4}-4tx_{2}y^{2}y_{3}z_{3}$ \\
$\hspace{0.55cm}-3tx_{2}yy_{2}y_{4}z+4tx_{2}yy_{3}^{2}z-3tx_{3}y^{2}y_{1}z_{4}+4tx_{3}y^{2}y_{3}z_{2}+tx_{3}y^{2}y_{4}z_{1}z$ \\
$\hspace{0.55cm}+2tx_{3}yy_{1}y_{4}z-4tx_{3}yy_{2}y_{3}+3tx_{4}y^{2}y_{1}z_{3}-3tx_{4}y^{2}y_{2}z_{2}-tx_{4}y^{2}y_{3}z_{1}$ \\
$\hspace{0.55cm}-2tx_{4}yy_{1}y_{3}z+3tx_{4}yy_{2}^{2}z-3x^{2}yy_{2}z_{2}z_{4}+4x^{2}yy_{2}z_{3}^{2}+2x^{2}yy_{3}z_{1}z_{4}$ \\
$\hspace{0.55cm}-4x^{2}yy_{3}z_{2}z_{3}-2x^{2}yy_{4}z_{1}z_{3}+3x^{2}yy_{4}z_{2}^{2}-2xx_{1}yy_{3}zz_{4}+2xx_{1}yy_{4}zz_{3}$ \\
$\hspace{0.55cm}+3xx_{2}y^{2}z_{2}z_{4}-4xx_{2}y^{2}z_{3}^{2}+3xx_{2}yy_{2}zz_{4}+4xx_{2}yy_{3}zz_{3}-6xx_{2}yy_{4}zz_{2}$ \\
$\hspace{0.55cm}-2xx_{3}y^{2}z_{1}z_{4}+4xx_{3}y^{2}z_{2}z_{3}-8xx_{3}yy_{2}zz_{3}+4xx_{3}yy_{3}zz_{2}+2xx_{3}yy_{4}zz_{1}$ \\
$\hspace{0.55cm}+2xx_{4}y^{2}z_{1}z_{3}-3xx_{4}y^{2}z_{2}^{2}+3xx_{4}yy_{2}zz_{2}-2xx_{4}yy_{3}zz_{1}+2x_{1}x_{3}y^{2}zz_{4}$ \\
$\hspace{0.55cm}-2x_{1}x_{3}yy_{4}z^{2}-2x_{1}x_{4}y^{2}zz_{3}+2x_{1}x_{4}yy_{3}z^{2}-3x_{2}^{2}y^{2}zz_{4}+3x_{2}^{2}yy_{4}z^{2}$ \\
$\hspace{0.55cm}+4x_{2}x_{3}y^{2}zz_{3}-4x_{2}x_{3}yy_{3}z^{2}+3x_{2}x_{4}y^{2}zz_{2}-3x_{2}x_{4}yy_{2}z^{2}-4x_{3}^{2}y^{2}zz_{2}$ \\
$\hspace{0.55cm}+4x_{3}^{2}yy_{2}z^{2}-xy^{2}y_{3}z_{4}+xy^{2}y_{4}z_{3}+x_{3}y^{3}z_{4}-x_{3}y^{2}y_{4}z-x_{4}y^{3}z_{3}+x_{4}y^{2}y_{3}z)$. \\ \\ \hline
\end{longtable}

\bigskip

\textsc{Acknowledgement.}
JOJ thanks Fredrik Andreassen, Eivind Schneider and Henrik Winther for useful discussions.
BK thanks Peter Olver and Niky Kamran for helpful correspondence.
The publication charges for this article have been funded by a grant from the publication fund of UiT The Arctic University of Norway.


\begin{thebibliography}{WW}
 \footnotesize

\bibitem{AK}
F.\ Andreassen, B.\ Kruglikov,
{\it Joint Invariants of Linear Symplectic Actions}, arXiv:2010.09899 (2020).

\bibitem{CW}
S.\,S.\ Chern, H.\,C.\ Wang, {\it Differential geometry in symplectic space. I},
Sci.\ Rep.\ Nat.\ Tsing Hua Univ.\ {\bf 4}, 453-477 (1947).

\bibitem{D}
V.\ Deconchy, {\it Hypersurfaces in symplectic affine geometry},
Differ.\ Geom.\ Appl.\ {\bf 17}, 1-13 (2002).

\bibitem{DFKN}
B.\ Doubrov, E.\ Ferapontov, B.\ Kruglikov, V.\ Novikov,
{\it On the integrability in Grassmann geometries: integrable systems associated with fourfolds Gr(3,5)},
Proc.\ London Math.\ Soc.\ {\bf 116}, 1269-1300 (2018).

\bibitem{JJ}
J.\,O.\ Jensen, {\it Differential Invariants of Symplectic and
Contact Lie alebra Actions}, Master Thesis in Mathematics,
hdl.handle.net/10037/19002,
UiT the Arctic University of Norway, June 2020.

\bibitem{KOT}
N.\ Kamran, P.\,J.\ Olver, K.\ Tenenblat, {\it Local symplectic invariants for curves},
Commun.\ Contemp.\ Math.\ {\bf 11}, 165-183 (2009).

\bibitem{Kr}
B.\ Kruglikov, {\it Poincar\'e function for moduli of differential-geometric structures},
Moscow Math.\ Journ.\ {\bf 19}, no.\ 4, 761-788 (2019).

\bibitem{KL0}
B.\ Kruglikov, V.\ Lychagin, {\it Geometry of Differential
equations\/}, in: Handbook of Global Analysis, Ed. D.Krupka, D.Saunders, Elsevier, 725-772 (2008).

\bibitem{KL1}
B.\ Kruglikov, V.\ Lychagin, {\it Differential invariants of the motion group actions},
in: Variations, Geometry and Physics, NOVA Scientific Publ.\ 237-251 (2009).

\bibitem{KL}
B. Kruglikov, V. Lychagin, {\it Global Lie-Tresse theorem}, Selecta Math.
{\bf 22}, 1357-1411 (2016).

\bibitem{MK}
B.\ McKay, {\it Lagrangian submanifolds in affine symplectic geometry},
Differ.\ Geom.\ Appl.\ {\bf 24} (6), 670-689 (2006).

\bibitem{MH}
E.\ Musso, E.\ Hubert, {\it Lagrangian curves in a 4-dimensional affine symplectic space},
Acta Appl.\ Math.\ {\bf 134}, 133-160 (2014).

\bibitem{MN}
E.\ Musso, L.\ Nicolodi, {\it Symplectic applicability of Lagrangian surfaces},
SIGMA Symmetry Integrability Geom.\ Methods Appl.\ {\bf 5}, Paper 067, 18 pp (2009).

\bibitem{Ol}
P.\,J.\ Olver, {\it Equivalence, Invariants and Symmetry}, Cambridge University Press (1995).

\bibitem{Ov}
L.\,V.\ Ovsiannikov, {\it Group Analysis of Differential Equations}, Academic Press (1982).

\bibitem{T}
T.Y. Thomas, {\it The Differential Invariants of Generalized Spaces},
Cambridge University Press, Cambridge (1934).

\bibitem{V}
F.\ Valiquette, {\it Geometric affine symplectic curve flows in $\R^4$},
Differ.\ Geom.\ Appl.\ {\bf 30} (6), 631-641 (2012).

\bibitem{X}
X.\ Xu, {\it Differential invariants of classical groups}, Duke Math.\ Journ.\ {\bf 94}, no.\ 3, 543-572 (1998).

\end{thebibliography}


\end{document}